\documentclass{amsart}
\usepackage{amsmath,amsthm,bbm}
\usepackage{amsfonts,amssymb}
\usepackage{epsfig}

 \newtheorem{theorem}{Theorem}[section]
 \newtheorem{lemma}[theorem]{Lemma}
 
 \newtheorem{proposition}[theorem]{Proposition}
 \newtheorem{definition}[theorem]{Definition}

 \newcounter{figures}[section]

\newcommand{\beq}{\begin{equation}}
\newcommand{\eeq}{\end{equation}}

\newcommand\norm[1]{\|#1\|}
\newcommand\nnorm[1]{\Big\|#1\Big\|}

\newcommand{\R}{\mathbb{R}}
\newcommand{\Rd}{\R^d}

\newcommand{\X}{\mathcal{X}}

  \def\cV{\mathcal{V}}

  \def\cN{\mathcal{N}}

  \def\cD{\mathcal{D}}
  
  \def\cO{\mathcal{O}}
  \def\cR{{\mathcal R}}
  
  \def\CD{\mathcal{D}}

\def\ip#1{{\langle #1 \rangle}}

\def\ONE{{\mathbbm 1}}
\def\tONE{{\tilde{\mathbbm 1}}}

\def\supp{\rm{supp}\, }

\newcommand{\cM}{{\mathcal M}}

\newcommand{\wh}{\widehat}
\newcommand{\ha}{\widehat a}
\newcommand{\hb}{\widehat b}

\newcommand\F[4]{F^{#1#2}_{#3#4}}
\newcommand\B[4]{B^{#1#2}_{#3#4}}

\def\Fsrpq{{F_{pq}^{s\r}}}
\def\Fspq{{F_{pq}^{s}}}
\def\fsrpq{{f_{pq}^{s\r}}}

\def\Bsrpq{{B_{pq}^{s\r}}}
\def\Bspq{{B_{pq}^{s}}}
\def\bsrpq{{b_{pq}^{s\r}}}
\def\Bstau{{B_\tau^{s}}}

\newcommand{\eps}{{\varepsilon}}

\def\cX{{\mathcal X}}

\def\wt{\widetilde}

\def\RR{{\mathbb R}}

\def\proj{\operatorname{Proj}}

\def\Ker{{L}}
\def\Lp{{p}}
\def\L2{{2}}
\def\Lpp{{L_p(\W)}}

\def\Linfty{{\infty}}
\def\bC{{\mathbb C}}
\def\CC{{\mathbb C}}
\def\ph{{\varphi}}

\def\V{V}

\def\W{{w_\mu}}
\def\WW{{W_\mu}}
\def\PP{{\mathsf {P}}}
\def\r{{\rho}}

\begin{document}

\title{Decomposition of weighted Triebel-Lizorkin and Besov spaces on the ball}

\author{George Kyriazis}
\address{Department of Mathematics and Statistics,
University of Cyprus, 1678 Nicosia, Cyprus}
\email{kyriazis@ucy.ac.cy}

\author{Pencho Petrushev}
\address{Department of Mathematics,
University of South Carolina, Columbia, SC 29208, USA}
\email{pencho@math.sc.edu}

\author{Yuan Xu}
\address{Department of Mathematics,
University of Oregon, Eugene, Oregon, 97403, USA}
\email{yuan@math.uoregon.edu}
\thanks{The third author has been partially supported by NSF Grant DMS-0604056.}

\thanks{Address for manuscript correspondence:
George Kyriazis, Department of Mathematics and Statistics,
University of Cyprus, P.O. Box 20537, 1678 Nicosia, Cyprus E-mail:
kyriazis@ucy.ac.cy }




\begin{abstract}
Weighted Triebel-Lizorkin and Besov spaces
on the unit ball $B^d$ in $\Rd$ with weights $\W(x)= (1-|x|^2)^{\mu-1/2}$, $\mu \ge 0$,
are introduced and explored.
A decomposition scheme is developed in terms of  almost exponentially localized
polynomial elements (needlets) $\{\varphi_\xi\}$, $\{\psi_\xi\}$
and it is shown that the membership of a distribution to
the weighted Triebel-Lizorkin or Besov spaces can be determined by the
size of the needlet coefficients $\{\ip{f,\varphi_\xi}\}$ in
appropriate sequence spaces.
\end{abstract}

\maketitle

\pagestyle{myheadings}
\thispagestyle{plain}
\markboth{GEORGE KYRIAZIS, PENCHO PETRUSHEV, AND YUAN XU}
         {WEIGHTED TRIEBEL-LIZORKIN AND BESOV SPACES ON THE BALL}

\section{Introduction}\label{introduction}
\setcounter{equation}{0}

Localized bases and frames allow to decompose functions and distributions
in terms of building blocks of simple nature and have numerous advantages
over other means of representation.
In particular, they enable one to encode smoothness and other norms
in terms of the coefficients of the decompositions.
Meyer's wavelets \cite{Meyer} and the $\varphi$-transform of
Frazier and Jawerth \cite{FJ1, FJ2, FJW} provide such building
blocks for decomposition of Triebel-Lizorkin and Besov spaces
in the classical case on $\R^d$.

The aim of this article is to develop similar tools
for decomposition of weighted Triebel-Lizorkin and Besov spaces
on the unit ball $B^d$ in $\RR^d$ ($d>1$) with weights
$$
\W(x):=(1-|x|^2)^{\mu-1/2}, \quad \mu \ge 0,
$$
were $|x|$ is the Euclidean norm of $x\in B^d$.
These include $L_p(B^d,\W)$, the Hardy spaces $H_p(B^d,\W)$, and
weighted Sobolev spaces.
For our purposes we develop localized frames which can be viewed as an analogue
of the $\varphi$-transform of Frazier and Jawerth on $B^d$.

For the construction of our frame elements we shall use
orthogonal polynomials in the weighted space $L_2(\W):= L_2(B^d, \W)$.
Denote by $\Pi_n$ the space of all algebraic polynomials of degree $n$
in $d$ variables and by $\V_n$ the subspace of all polynomials of degree $n$
which are orthogonal to lower degree polynomials in $L_2(\W )$.
These are eigenspaces of the differential operator
\begin{equation}\label{def-diff-oper}
D_\mu:= -\Delta +\langle x, \nabla \rangle^2
+(2\mu+d-1)\langle x, \nabla \rangle.
\end{equation}
More precisely (see e.g. \cite{DXu}),
\begin{equation}\label{eigen-func}
D_\mu P=n(n+d+2\mu -1)P
\quad\mbox{for } P\in\V_n.
\end{equation}
We have the orthogonal polynomial decomposition
\begin{equation}\label{L2:decomp}
L_2(\W) = \bigoplus_{n=0}^\infty \V_n, \qquad
\V_n \subset \Pi_n.
\end{equation}
Note that $\dim \V_n = \binom{n+d-1}{n}\sim n^{d-1}$.
As is shown in \cite{X99} the orthogonal projector
$\proj_n: L_2(\W) \mapsto\V_n$
can be written as
\begin{equation}\label{def-Pn}
(\proj_n f)(x) = \int_{B^d} f(y)\PP_n(x,y) \W(y) dy,
\end{equation}
where, for $\mu > 0$, the kernel $\PP_n(x,y)$ has the representation
\begin{align} \label{compactPn}
\PP_n(x,y)& = b_d^\mu b_1^{\mu-\frac{1}{2}}\frac{n+\lambda}{\lambda} \\
&\times  \int_{-1}^1 C_n^\lambda \left(\langle x,y\rangle +
u\sqrt{1-|x|^2} \sqrt{1-|y|^2}\right) (1-u^2)^{\mu-1}du.  \notag
\end{align}
Here $\langle x, y \rangle$ is the Euclidean inner product in $\R^d$,
$C_n^\lambda$ is the $n$-th degree Gegenbauer polynomial,
\begin{equation}\label{def-lambda}
   \lambda = \mu + \frac{d-1}{2},
\end{equation}
and the constants
$b_d^\mu$, $ b_1^{\mu-\frac{1}{2}}$ are defined by
$(b_d^\gamma)^{-1} := \int_{B^d} (1-|x|^2)^{\gamma-1/2} dx$.
For a~representation of $\PP_n(x, y)$ in the limiting case $\mu=0$,
see (4.2) in \cite{PX2}.

Evidently,
\begin{equation}\label{def:Kn}
K_n(x,y):= \sum_{j=0}^n \PP_j(x,y)
\end{equation}
is the kernel of the orthogonal
projector of $L_2(\W)$ onto the space
$\bigoplus_{\nu=0}^n \V_\nu$.

A key role in this study will play the fact (established in \cite{PX2})
that if the coefficients on the right in (\ref{def:Kn})
are ``smoothed out" by sampling a compactly supported $C^\infty$ function,
then the resulting kernel has nearly exponential localization
around the main diagonal $y=x$ in $B^d\times B^d$.
More precisely, let
\begin{equation} \label{def-Ln}
L_n(x,y):= \sum_{j=0}^\infty \ha\Big(\frac{j}{n}\Big) \PP_j(x,y),
\end{equation}
where the ``smoothing" function $\ha$ is admissible in the sense of the following definition:


\begin{definition} \label{defn:admissible}
A function $\ha \in C^\infty[0, \infty)$
is called admissible of type

$(a)$ if
$\supp \ha \subset [0,2]$ and $\ha(t)=1$ on $[0, 1]$,
and of type

$(b)$ if $\supp \ha \subset [1/2,2]$.
\end{definition}
\noindent
We introduce the distance
\begin{equation} \label{eq:distant}
d(x,y):= \arccos \left
\{ \langle x,y\rangle + \sqrt{1-|x|^2}\sqrt{1-|y|^2} \right \}
\quad \mbox{on $B^d$}
\end{equation}
and set
\begin{equation} \label{def-WW}
\WW(n;x) 
:= \left(\sqrt{1-|x|^2} + n^{-1}\right)^{2\mu}, \quad x \in B^d.
\end{equation}
One of our main results in \cite[Theorem 4.2]{PX2} asserts that for any $k >0$
there exists a constant $c_k>0$ depending only on $k$, $d$, $\mu$, and $\wh a$
such that
\begin{equation} \label{eq:main-estI}
|L_n(x,y)| \le c_k \frac{n^d }{\sqrt{\WW(n;x)} \sqrt{\WW(n;y)}
(1 + n\,d(x,y))^k},
\quad x, y\in B^d.
\end{equation}

The kernels $L_n$ are our main ingredient in constructing
{\it analysis} and {\it synthesis} needlet systems
$\{\ph_\xi\}_{\xi\in\cX}$ and $\{\psi_\xi\}_{\xi\in\cX}$ here,
indexed by a multilevel set $\cX=\cup_{j=0}^\infty \cX_j$ (\S\ref{def-needlets}).
This is a pair of dual frames whose elements have nearly exponential localization on $B^d$
and provide representation of every distribution $f$ on $B^d$:
\begin{equation}\label{eq:decomposition}
f=\sum_{\xi \in \cX}\langle f, \ph_\xi \rangle \psi_\xi.
\end{equation}
The superb localization of the frame elements prompted us to term them {\em needlets}.

Our main interest lies with distributions in the weighted Triebel-Lizorkin ($F$-spaces) and
 Besov spaces ($B$-spaces) on $B^d$.
These spaces are naturally defined via spectral decompositions
(see \cite{Pee, T1} for the general idea).
To be specific, let
$$
\Phi_0(x, y) := 1
\quad\mbox{and}\quad
\Phi_j(x, y) := \sum_{\nu=0}^\infty \ha
\Big(\frac{\nu}{2^{j-1}}\Big)\PP_\nu(x,y), \quad j\ge 1,
$$
where $\PP(\cdot, \cdot)$ is from \eqref{compactPn} and
$\ha$ is admissible of type (b) (see Definition~\ref{defn:admissible})
such that $|\ha|>0$ on $[3/5, 5/3]$.

The $F$-space $\Fsrpq$ with
$s, \r \in \R$, $0<p<\infty$, $0<q\le\infty$,
is defined (\S\ref{Tri-Liz})
as the space of all distributions $f$ on $B^d$ such that
\begin{equation}\label{Tri-Liz-norm0}
\norm{f}_{\Fsrpq}:=\nnorm{\biggl(\sum_{j=0}^{\infty}
(2^{sj}\WW(2^j; \cdot)^{-\r/d}|\Phi_j\ast f(\cdot)|)^q\biggr)^{1/q}}_{\Lpp} <\infty,
\end{equation}
where $\Phi_j*f(x):= \langle f, \overline{\Phi(x, \cdot)}\rangle$
(see Definition~\ref{Def-conv}).

The corresponding scales of weighted Besov spaces $\Bsrpq$
with $s, \r \in \R$, $0<p, q\le\infty$, are
defined (\S\ref{Besov}) via the (quasi-)norms
\begin{equation}\label{Besov-norm0}
\|f\|_{\Bsrpq} :=
\Big(\sum_{j=0}^\infty \Big(2^{s j}
\|\WW(2^j; \cdot)^{-\r/d}\Phi_j*f(\cdot)\|_{\Lpp}\Big)^q\Big)^{1/q}.
\end{equation}

Unlike in the classical case on $\R^d$, we have introduced an additional
parameter~$\r$, which allows considering different scales
of Triebel-Lizorkin and Besov spaces.
To~us most natural are the spaces
\begin{equation}\label{Tri-Liz-Besov}
\Fspq:=F_{pq}^{s s}\quad \mbox{and}\quad \Bspq:=B_{pq}^{s s},
\end{equation}
which scale (are embedded) correctly with respect to the smoothness parameter~$s$.
A ``classical" choice would be to consider the spaces $F_{pq}^{s 0}$ and $B_{pq}^{s 0}$,
where the weight $\WW(2^j;\cdot)$ is excluded from
\eqref{Tri-Liz-norm0}-\eqref{Besov-norm0}.
The introduction of the parameter $\r$ enables us to treat these
spaces simultaneously.

One of the main results of this paper is the characterization of  the $F$-spaces
in terms of the size of the needlet coefficients in the decomposition (\ref{eq:decomposition}),
namely,
$$
\norm{f}_{\Fsrpq}
\sim \nnorm{\biggl(\sum_{j=0}^\infty 2^{sjq}
\sum_{\xi\in \cX_j}|\ip{f, \varphi_\xi}|\WW(2^j;\xi)^{-\r/d}|\psi_\xi(\cdot)|^q\biggr)^{1/q}}_{\Lpp}.
$$
Similarly for the Besov spaces $\Bsrpq$ we have the characterization (\S\ref{Besov})
$$
\norm{f}_{\Bsrpq}
\sim \Big(\sum_{j=0}^\infty2^{sjq}
\Bigl[\sum_{\xi\in \cX_j}
\Big(\WW(2^j; \xi)^{-\r/d}\norm{\ip{f,\varphi_\xi}\psi_\xi}_{\Lpp}\Big)^p\Bigr]^{q/p}\Bigr)^{1/q}.
$$

Further, the weighted Besov spaces are applied to nonlinear $n$-term approximation from
needlets on $B^d$ (\S\ref{Nonlin-app}).


This is a follow-up paper of \cite{PX2}, where the localization \eqref{eq:main-estI}
is established and the construction and basic properties of a single system of needlets are given.
Our development here is a part of a broader undertaking for needlet characterization of
Triebel-Lizorkin and Besov spaces on nonclassical domains,
including the multidimensional unit sphere \cite{NPW1, NPW2}, ball,
and cube (interval \cite{KPX, PX1}) with weights.
The~results in this paper generalize the results in the univariate case from \cite{KPX}
(with $\alpha=\beta$), where needlet characterizations of $F$- and $B$-spaces on
the interval are obtained.
%

The organization of the paper is the following:
In \S\ref{preliminaries} the needed results from \cite{PX2} and some background material
are given, including localized polynomial kernels, the maximal operator,
distributions on $B^d$, and cubature formula on $B^d$.
The definition and some basic properties of needlets are given in \S\ref{def-needlets}.
In \S\ref{Tri-Liz} the weighted Triebel-Lizorlin space on $B^d$ are introduced
and characterized via needlets,
while the weighted Besov spaces are explored in \S\ref{Besov}.
In \S\ref{Nonlin-app} Besov spaces are applied to nonlinear $n$-term approximation
from needlets.
Section~\ref{proofs} contains the proofs of various lemmas from previous sections.

Throughout the paper we use the following notation:
$$
\|f\|_\Lp:=\Big(\int_{B^d}|f(x)|^p\W(x)dx\Big)^{1/p},
\quad 0<p<\infty, \quad 
\|f\|_\Linfty:={\rm ess \; sup}_{x\in B^d}|f(x)|.
$$
For a measurable set $E\subset B^d$,
$|E|$ denotes the Lebesgue measure of $E$,
$m(E):= \int_E\W(x)dx$,
$\ONE_E$ is the characteristic function of $E$, and
$\tONE_E:=|m(E)|^{-1/2}\ONE_E$ is the $L_2(\W)$ normalized
characteristic function of $E$.
Positive constants are denoted by $c$, $c_1, c_*, \dots$ and they
may vary at every occurrence; $A\sim B$ means $c_1A\le B\le c_2 A$.

\section{Preliminaries}\label{preliminaries}
\setcounter{equation}{0}

\subsection{Localized polynomial kernels on \boldmath $B^d$}\label{Local-kernels}

The polynomial kernels $L_n(x,y)$ introduced in (\ref{def-Ln})
will be our main vehicle in developing needlet systems.
Here we give come additional properties of these kernels.

%
We have
\begin{equation} \label{eq:LpUB}
\|L_n(x, \cdot)\|_{\Lp}
\le c\Bigl(\frac{n^{d}}{\WW(n; x)}\Bigr)^{1-1/p},
\quad x\in B^d,
\quad 0<p\le\infty.
\end{equation}
This estimate is an immediate consequence of \eqref{eq:main-estI}
and the following lemma (see \cite[Lemma 4.6]{PX2}),
which will be instrumental in several proofs below.


\begin{lemma}\label{lem:instrumental}
If $\sigma > d/p + 2\mu|1/p-1/2|$, $\mu\ge 0$, $0 < p < \infty$, then
\begin{equation} \label{est-inst}
\int_{B^d} \frac{\W(y) dy}
{\WW(n;y)^{p/2} (1+n d(x,y))^{\sigma p}}
\le c\,n^{-d}\WW(n;x)^{1-p/2}.
\end{equation}
\end{lemma}

We now establish a matching lower bound estimate. 


\begin{theorem}\label{thm:est-Lp-norm}
Let $\ha$ be admissible and $|\ha(t)| \ge c_* > 0$ for $t \in [3/5,5/3]$.
Then for $0<p\le \infty$ and $n\ge 2$
\begin{equation} \label{eq:est-Lp-norm}
\|L_n(x, \cdot)\|_{\Lp}
\ge c\Bigl(\frac{n^{d}}{\WW(n; x)}\Bigr)^{1-1/p},
\quad x \in B^d.
\end{equation}
Here the constant $c>0$ depends only on $d$, $\mu$, $p$, and $c_*$.
\end{theorem}

The proof of this theorem is given in \S\ref{proofsA}.

The kernels $L_n(x,y)$
are in a sense Lip1 functions in both variables with respect
to the distance $d(\cdot,\cdot)$ from \eqref{eq:distant}:
Let $\xi, y\in B^d$ and $c^*>0$, $n\ge 1$.
Then for all $x, z\in B_\xi(c^*n^{-1})$ and an arbitrary $k$, we have
\begin{equation} \label{Lip}
|L_n(x,y)-L_n(\xi,y)|
\le c_k \frac{n^{d+1}d(x, \xi)} {\sqrt{\WW(n; y)}\sqrt{\WW(n; z)}(1+nd(y,z))^k},
\end{equation}
where $c_k$ depends only on $k$, $\mu$, $d$, $\ha$, and $c^*$ (see
\cite[Proposition 4.7]{PX2}).

\medskip
We shall also need the following inequality from \cite[Lemma 4.1]{PX2}:
\begin{equation} \label{norm-dist2}
\Big|\sqrt{1-|x|^2}-\sqrt{1-|y|^2}\Big| \le \sqrt{2}\, d(x, y),
\quad x,y\in B^d,
\end{equation}
which yields
\begin{equation}\label{eq:useful}
\WW(n;x)\le 2^\mu\WW(n;y)(1+nd(x,y))^{2\mu},
\quad x,y\in B^d.
\end{equation}


\subsection{Reproducing polynomial kernels and applications}\label{kernels}

To simplify our notation we introduce the following
``convolution": For functions $\Phi: B^d\times B^d \to \bC$ and
$f: B^d \to \bC$, we write
\begin{equation}\label{convolution}
\Phi*f(x) := \int_{B^d} \Phi(x, y)f(y)\W(y)\,dy.
\end{equation}
We denote by $E_n(f)_p$ the best approximation
of $f \in \Lpp$ from $\Pi_n$, i.e.
\begin{equation}\label{def-En}
E_n(f)_p := \inf_{g \in \Pi_n}\|f-g\|_\Lp.
\end{equation}


\begin{lemma}\label{lem:Ker-n}
Let $L_n$ be the kernel from $(\ref{def-Ln})$,
with $\ha$ admissible of type $(a)$. Then

$(i)$ $\Ker_n*g =g$ for $g \in \Pi_n$, i.e.
$\Ker_n$ is a reproducing kernel for
$\Pi_n$, and

$(ii)$ for any $f \in \Lpp$, $1\le p \le \infty$, we have
$\Ker_n*f \in \Pi_{2n}$,
\begin{equation}\label{Ker-n3}
\|\Ker_n*f\|_\Lp\le c \|f\|_\Lp, \quad\mbox{and}\quad
\|f-\Ker_n*f\|_\Lp \le cE_n(f)_p.
\end{equation}
\end{lemma}

This lemma follows readily by the definition of $L_n$
(see also Definition~\ref{defn:admissible}) and (\ref{eq:LpUB})
(see \cite[Proposition 4.8]{PX2}).

Lemma~\ref{lem:Ker-n} (i) and (\ref{eq:LpUB}) are instrumental
in relating weighted norms of polynomials.


\begin{proposition}\label{Nikolski}
For $0 < q \le p \le \infty$ and $g \in \Pi_n$, $n\ge 1$,
\begin{equation}\label{norm-relation}
\|g\|_p \le cn^{(d+2\mu)(1/q-1/p)}\|g\|_q,
\end{equation}
and for any $\gamma\in\RR$
\begin{equation}\label{norm-relation2}
\|\WW(n;\cdot)^\gamma g(\cdot)\|_p \le
cn^{d(1/q-1/p)}\|\WW(n;\cdot)^{\gamma+1/p-1/q} g(\cdot)\|_q.
\end{equation}
\end{proposition}

The proof of this proposition is quite similar to the proof of Proposition 2.6
from~\cite{KPX}; for completeness it is given in \S\ref{proofsA}.

\subsection{Maximal operator}\label{Max-iequal}

We denote by $B_\xi(r)$ the ball centered at $\xi\in B^d$ of
radius $r>0$ with respect to the distance $d(\cdot,\cdot)$ on
$B^d$, i.e.
\begin{equation}\label{def-ball}
B_\xi(r)= \{x\in B^d: d(x, \xi)< r\}.
\end{equation}
It is straightforward to show that (see \cite[Lemma 5.3]{PX2})
\begin{equation}\label{ball-measure}
|B_\xi(r)|\sim r^d\sqrt{1-|\xi|^2}
\end{equation}
and
\begin{equation}\label{ball-size}
m(B_\xi(r)):=\int_{B_\xi(r)}\W(x)\, dx
\sim r^d(r+\sqrt{1-|\xi|^2})^{2\mu}
\sim r^d(r+ d(\xi,\partial B^d))^{2\mu},
\end{equation}
where $\partial B^d$ is the boundary of $B^d$,
i.e. the unit sphere in $\R^d$.

The maximal operator $\cM_t$ ($t>0$) is defined by
\begin{equation}\label{def: max-op}
\cM_tf(x):=\sup_{B\ni x}\left(\frac1{m(B)}\int_B|f(y)|^t\W(y)\, dy\right)^{1/t},
\quad x\in B^d,
\end{equation}
where the sup is over all balls  $B\subset B^d$ (with respect to $d(\cdot, \cdot)$)
containing~$x$.

It follows by \eqref{ball-size} that the measure
$m(E):=\int_E\W(x)\, dx$
is a doubling measure on $B^d$,
i.e. for $\xi\in B^d$ and $0<r\le \pi$
\begin{equation}\label{eq:doub-vol}
m(B_\xi(2r))\le c m(B_\xi(r)).
\end{equation}
Consequently, the general theory of maximal operators applies
and the Fefferman-Stein vector-valued maximal inequality is valid (see \cite{Stein}):
If $0<p<\infty, 0<q\le\infty$, and
$0<t<\min\{p,q\}$ then for any sequence of functions
$\{f_\nu\}_\nu$ on $B^d$
\begin{equation}\label{max-ineq}
\nnorm{\Bigl(\sum_{\nu=1}^\infty|\cM_tf_\nu(\cdot)|^q\Bigr)^{1/q}}_\Lp
\le \nnorm{\Bigl(\sum_{\nu=1}^\infty|
f_\nu(\cdot)|^q\Bigr)^{1/q}}_\Lp.
\end{equation}

We need to estimate $\cM_t \ONE_{B}$ for an arbitrary ball $B\subset B^d$.

\begin{lemma}\label{lem:B-maximal}
Let $\xi \in B^d$ and $0<r\le \pi$.
Then for $x\in B^d$
\begin{equation}\label{B-max1}
(\cM_t \ONE_{B_\xi(r)})(x)
\sim \Big(1+\frac{d(\xi , x)}{r}\Big)^{-d/t}
     \Big(1+\frac{d(\xi , x)}{r+d(\xi , \partial B^d)}\Big)^{-2\mu/t}
\end{equation}
and hence
\begin{equation}\label{B-max3}
c'\Big(1+\frac{d(\xi , x)}{r}\Big)^{-(2\mu+d)/t}
\le (\cM_t \ONE_{B_\xi (r)})(x)
\le c\Big(1+\frac{d(\xi , x)}{r}\Big)^{-d/t}.
\end{equation}
Here the constants depend only on $d$, $\mu$, and $t$.
\end{lemma}

\begin{proof}It is easy to see that
$$
(\cM_t \ONE_{B_\xi (r)})(x)
=\sup_{B\ni x}\left(\frac{m(B\cap B_{\xi}(r))}{m(B)}\right)^{1/t},
\quad x\in B^d,
$$
where the sup is taken over all the balls $B\subset B^d$
(with respect to $d(\cdot, \cdot)$) containing~$x$.
This immediately leads to $(\cM_t \ONE_{B_\xi (r)})(x)\sim 1$ if $d(x, \xi)\le 2r$
and hence (\ref{B-max1}) holds in this case.

Suppose $d(\xi,x)>2r$. Then evidently
$$
(\cM_t \ONE_{B_\xi (r)})(x)^{1/t} \ge \left(\frac{m(B_{\xi}(r))}{m(B_{\xi}(d(x,\xi)))}\right)^{1/t}.
$$
For the other direction, suppose $B_z(r^*)\subset B^d$ is the smallest ball such that
$x\in \overline{B_z(r^*)}$ and $\overline{B_z(r^*)}\cap \overline{B_\xi(r)}\ne \emptyset$.
A simple application of the triangle inequality shows that
$B_\xi(d(\xi,x))\subset B_z(5r^*)$.
Thus using (\ref{eq:doub-vol})
$$
(\cM_t \ONE_{B_\xi (r)})(x)\le \left(\frac{m(  B_{\xi}(r))}{m(B_z(r^*))}\right)^{1/t}
\le c \left(\frac{m(B_{\xi}(r))}{m(B_{\xi}(d(x,\xi))}\right)^{1/t}.
$$
Therefore, using \eqref{ball-size}
$$
(\cM_t \ONE_{B_\xi (r)})(x)
\sim \left(\frac{m(B_{\xi}(r))}{m(B_{\xi}(d(x,\xi))}\right)^{1/t}
\sim \left(\frac{r^d(r+d(\xi, \partial B^d))^{2\mu}}
                {d(x,\xi)^d(d(x,\xi)+d(\xi, \partial B^d))^{2\mu} }\right)^{1/t},
$$
which implies (\ref{B-max1}) since $d(\xi,x)>2r$. Estimate
(\ref{B-max3}) is immediate from  (\ref{B-max1}).
\end{proof}


\subsection{Distributions on \boldmath $B^d$}\label{distributions}

To define distributions on $B^d$ we shall use as
test functions the set $\cD:=C^{\infty}(B^d)$ of all infinitely continuously differentiable
complex valued functions on $B^d$ such that
\begin{equation}\label{norms}
\|\phi\|_{W_\infty^k}
:=\sum_{|\alpha|\le k}\|\partial^\alpha \phi\|_\infty<\infty
\quad\mbox{for } k=0, 1, \dots.
\end{equation}
We assume that the topology in $\cD$ is defined by these norms.

Evidently all polynomials belong to $\cD$.   More importantly, the space $\cD$
of test functions $\phi$ can be completely characterized by their orthogonal
polynomial expansions. Denote
\beq\label{D-norms}
\cN_k(\phi):=\sup_{n\ge 0}\, (n+1)^k \|\proj_n \phi\|_2.
\eeq


\begin{lemma}\label{lem:char-D}
$(a)$
$\phi\in\cD$ if and only if
$\|\proj_n \phi\|_2 =\cO(n^{-k})$ for all $k$.

$(b)$
For each $\phi\in \cD$,
$
\phi=\sum_{n=0}^\infty \proj_n \phi,
$
where the convergence is in the topology of $\cD$.

$(c)$
The topology in $\cD$ can be equivalently defined by the norms
$\cN_k(\cdot)$, $k=0, 1, \dots$.
\end{lemma}


\begin{proof}Let $\phi\in\cD$. Assume that $Q_{n-1}\in\Pi_{n-1}$ ($n\ge
1$) is the polynomial of best $L_2(\W)$-approximation to $\phi$,
i.e. $\|\phi-Q_{n-1}\|_2=E_{n-1}(\phi)_2$. Since $\PP_n(x, \cdot)$
is orthogonal to $\Pi_{n-1}$,
$$
|\proj_n \phi (x)| = |\langle \phi, \PP_n(x,\cdot) \rangle|
=|\langle \phi-Q_{n-1}, \PP_n(x,\cdot)\rangle|
\le E_{n-1}(\phi)_2 \PP_n(x,x)^{1/2}.
$$
By the Jackson type estimate from \cite{X05}, for any $k\ge 1$,
$$
E_n(\phi)_2
\le c_k n^{-2k} \|D_\mu^k\phi\|_{2}
\le c n^{-2k} \|D_\mu^k\phi\|_\infty
\le c n^{-2k}\sum_{|\alpha|\le 2k}\|\partial^\alpha \phi\|_\infty
= c n^{-2k}\|\phi\|_{W_\infty^{2k}}.
$$
Here $D_\mu$ is the differential operator from \eqref{def-diff-oper}.
It is easy to see that
$$
\|\PP_n(x,x)^{1/2}\|_2^2 = \binom{n+d-1}{n} \sim n^{d-1}.
$$
All of the above leads to
$$
\|\proj_n \phi\|_2  \le c_k n^{-2 k+ (d-1)/2} \|\phi\|_{W_\infty^{2k}},
\quad n\ge 1,
\quad \mbox{for any } k\ge 1.
$$
Therefore, for any $m\ge 0$
$$
\cN_m(\phi) \le c\|\phi\|_{W_\infty^{2k}}
\quad\mbox{if \; $k\ge m/2+(d-1)/4$.}
$$

In the other direction, by Markov's inequality (see \cite{Kellogg}) and \eqref{norm-relation},
it follows that
$$
\|\partial^\alpha \proj_n \phi \|_\infty
\le n^{2|\alpha|}\|\proj_n \phi\|_\infty
\le cn^{2|\alpha|+d/2+\mu}\|\proj_n \phi\|_2.
$$
Consequently, if $\|\proj_n \phi\|_2 =\cO(n^{-k})$ for all $k$,
then
$
\partial^\alpha \phi=\sum_{n=0}^\infty  \partial^\alpha \proj_n \phi
$
for all multi-indices $\alpha$ with the series converging uniformly and
$$
\|\phi\|_{W_\infty^k}
\le c\sum_{|\alpha|\le k}\sum_{n=0}^\infty n^{2 |\alpha|+d/2+\mu} \|\proj_n \phi \|_2
\le c \cN_{m}(\phi),
\quad m\ge 2k+d/2+\mu+2.
$$
This completes the proof of the lemma.
\end{proof}

\medskip

The space $\cD':=\cD'(B^d)$ of distributions on $B^d$ is defined
as the set of all continuous linear functionals on $\cD$.
The pairing of $f\in \cD'$ and $\phi\in\cD$ will be denoted by
$\langle f, \phi \rangle:= f(\overline{\phi})$, which will be shown
to be consistent with the inner product
$
\langle f, g \rangle :=\int_{B^d} f(x)\overline{g(x)} \W(x)dx
$
in $L_2(\W)$.

We now extend the definition of the nonstandard ``convolution" from
\eqref{convolution} to distributions.


\begin{definition}\label{Def-conv}
Let $f \in \CD'$ and assume that $\Phi: B^d\times B^d \mapsto \CC$ is such that
$\Phi(x, \cdot)\in \cD$ for all $x\in B^d$.
We define
$$
( \Phi \ast f)(x) := \langle f, \overline{\Phi(x,\cdot)}\rangle,
$$
where on the right $f$ acts on $\overline{\Phi(x, y)}$ as a function of $y$.
\end{definition}

For later use we next record some simple properties of this ``convolution".


\begin{lemma}\label{lem:prop-conv}
$(i)$ If $f\in \cD'$ and $\Phi(\cdot,\cdot)\in C^\infty(B^d\times B^d)$, then
$\Phi*f\in \cD$,
and in particular $\PP_n*f \in \V_n$.
We define $\proj_n f:=\PP_n*f$.

$(ii)$
If $f\in \cD'$ and $\Phi(\cdot,\cdot)\in C^\infty(B^d\times B^d)$, then
$$
\langle \Phi*f, \phi \rangle = \langle f, \overline{\Phi}*\phi\rangle,
\quad \phi\in\cD.
$$

$(iii)$
Let $\Phi(\cdot,\cdot), \Psi(\cdot, \cdot)\in C^\infty(B^d\times B^d)$,
and $\Phi(x, y)=\Phi(y, x)$ and $\Psi(x, y)=\Psi(y, x)$ for $x, y\in B^d$.
Then for any $f\in\cD'$ and $x\in B^d$
$$
\Psi*\overline{\Phi}*f(x)=\langle\Psi(x, \cdot), \Phi(\cdot, \cdot)\rangle*f.
$$
\end{lemma}

The proof of this lemma is standard and will be omitted.

We next give the representation of distributions from $\cD'$
in terms of orthogonal polynomials on $B^d$.


\begin{lemma}\label{lem:dec-D1}
$(a)$
A linear functional $f\in\cD'$
if and only if there exists $k\ge 0$ such that
\begin{equation}\label{D1}
|\langle f, \phi\rangle|\le c_k\cN_k(\phi)
\quad \mbox{for all } \quad \phi \in \cD,
\end{equation}
Hence, for $f \in \cD'$ there exits $k\ge 0$ such that
\begin{equation}\label{D2}
\|\proj_n f\|_2= \|\PP_n*f\|_2
\le c_k(n+1)^k, \quad n=0, 1, \dots.
\end{equation}

$(b)$
Every $f\in\cD'$ has the representation
$
f=\sum_{n=0}^\infty \proj_n f
$
in distributional sense, i.e.
\begin{equation}\label{D3}
\langle f, \phi\rangle
=\sum_{n=0}^\infty  \langle \proj_n f, \phi  \rangle
=\sum_{n=0}^\infty  \langle \proj_n f, \proj_n\phi  \rangle
\quad\mbox{for all}
\quad \phi\in\cD,
\end{equation}
where the series converges absolutely.
\end{lemma}

\noindent
\begin{proof}(a) Part (a) follows immediately by the fact that the
topology in $\cD$ can be defined by the norms $\cN_k(\cdot)$
defined in \eqref{D-norms}.

(b) Using Lemma~\ref{lem:char-D} (b) we get for $\phi\in\cD$,
$$
\langle f, \phi \rangle
=\lim_{N\to\infty}\Big\langle f, \sum_{n=0}^N \proj_n \phi \Big\rangle
=\lim_{N\to\infty}\sum_{n=0}^N  \langle f, \proj_n \phi\rangle
=\sum_{n=0}^\infty \langle \proj_n f, \proj_n \phi \rangle,
$$
where the last equality is justified by using (\ref{D2}) and the
rapid decay of $\|\proj_n \phi\|_2$.
\end{proof}

\subsection{Cubature formula and subdivision of \boldmath $B^d$}
\label{cubature}
For the construction of our building blocks (needlets)
we shall  utilize the positive cubature formula given in \cite{PX2}. This formula is based
on almost equally distributed knots on $B^d$ with respect to the distance $d(\cdot,\cdot)$.



\begin{definition}\label{def.eq-dist-p}
We say that a set $\cX_\eps \subset B^d$, along with an associated
partition $\cR_\eps$ of $B^d$ consisting of measurable subsets of $B^d$,
is a {\em set of almost uniformly $\eps$-distributed points} on $B^d$
if
\begin{enumerate}
\item[(i)]
$B^d = \bigcup_{R \in \cR_\eps} R$ and the sets in $\cR_\eps$
do not overlap
$(R_1^\circ \cap R_2^\circ = \emptyset$ if $R_1 \ne R_2)$.
\item[(ii)]
For each $R \in \cR_\eps$ there is a unique $\xi \in \cX_\eps$
such that
$
B_\xi(c^*\eps) \subset R \subset B_\xi(\eps).
$
\end{enumerate}
Hence
$\# \cX_\eps = \# \cR_\eps \le c^{**}\eps^{-d}$.
Here the constant $c^*>0$, depending only on $d$, is fixed but sufficiently small, so that the existence
of sets of almost uniformly $\eps$-distributed points on $B^d$
is guaranteed $($see the next lemma$)$.
\end{definition}


\begin{lemma}\label{lem:points} \cite{PX2}
For a sufficiently small constant $c^*>0$, depending only on $d$,
and an arbitrary $0<\eps \le \pi$ there exists
a set $\cX_\eps \subset B^d$  of almost uniformly $\eps$-distributed points on $B^d$,
where the associated partition $\cR_\eps$ of $B^d$
consists of projections of spherical simplices.
\end{lemma}

An important element in the construction of needlets will be
the cubature formula given in \cite[Corollary 5.10]{PX2}:


\begin{proposition}\label{prop:cubature}
There exists a constant $c^\diamond >0$ {\rm(}depending only on
$d${\rm)} and a sequence $\{\cX_j\}_{j=0}^\infty$ of almost
uniformly $\eps_j$-distributed points on $B^d$ with
$\eps_j:=c^\diamond 2^{-j}$, and there exist positive coefficients
$\{\lambda_\xi\}_{\xi \in \cX_j}$ such that the cubature formula
\begin{equation}\label{eq:cubature}
\int_{B^d} f(x)\W(x)\,dx \sim \sum_{\xi \in \cX_j} \lambda_\xi
f(\xi)
\end{equation}
is exact for all polynomials of degree $\le 2^{j+2}$. In addition,
\begin{equation}\label{eq:coeff}
\lambda_\xi \sim 2^{-jd}\WW(2^j;\xi) 
\sim m(B_\xi(2^{-j}))
\end{equation}
with constants of equivalence depending only on $\mu$ and $d$.
\end{proposition}

It follows from above that
\begin{equation}\label{eq:mR}
m(R_\xi)\sim 2^{-jd}\WW(2^j;\xi)\sim \lambda_\xi,
\quad \xi \in \cX_j,
\end{equation}
while
\begin{equation}\label{mes-R-xi}
|R_\xi|\sim 2^{-jd}(\sqrt{1-|\xi|^2}+2^{-j}),
\quad \xi \in \cX_j.
\end{equation}

\section{Localized building blocks (Needlets) on \boldmath $B^d$}\label{def-needlets}
\setcounter{equation}{0}

We utilize the ideas from \cite{NPW2, KPX} in constructing a pair of sequences
of ``analysis" and ``synthesis" needlets on $B^d$.
Let $\ha,\hb$ satisfy the conditions
\begin{align}
&\quad \ha,\hb\in C^\infty[0, \infty),
\quad \supp \ha,\hb \subset [1/2, 2], \label{ha-hb1}\\
&\quad  |\ha(t)|, |\hb(t)| >c>0, \quad \text{if  } t \in [3/5, 5/3],\label{ha-hb2}\\
&\quad \overline{\ha(t)}\, \hb(t) + \overline{\ha(2t)}\,\hb(2t)  =1,\label{ha-hb3}
\quad\text{if } t \in [1/2, 1].
\end{align}
Hence,
\beq\label{part-unity}
\sum_{\nu=0}^\infty \overline{\ha(2^{-\nu}t )}\, \hb(2^{-\nu}t)=1,
\quad t\in [1,\infty).
\eeq
It is easy to see that
if $\ha$ satisfies (\ref{ha-hb1})-(\ref{ha-hb2}), then
there exists $\hb$ satisfying (\ref{ha-hb1})-(\ref{ha-hb2})
such that (\ref{ha-hb3}) is valid (see e.g. \cite{FJ2}).

\smallskip

Let $\ha$, $\hb$ satisfy (\ref{ha-hb1})-(\ref{ha-hb3}).
We define
$\Phi_0(x, y)= \Psi_0(x, y):=1$,
\begin{equation}\label{def.Phi-j}
\Phi_j(x, y) := \sum_{\nu=0}^\infty
\ha\Big(\frac{\nu}{2^{j-1}}\Big) \PP_\nu(x,y), \quad j\ge 1,
\end{equation}
\begin{equation}\label{def-Psi-j}
\Psi_j(x, y) := \sum_{\nu=0}^\infty
\hb\Big(\frac{\nu}{2^{j-1}}\Big) \PP_\nu(x,y), \quad j\ge 1.
\end{equation}

Assume that $\cX_j$ is the set of knots and $\lambda_\xi$'s are the coefficients of
the  cubature formula \eqref{eq:cubature}.
We define the $j$th level {\em needlets} by
\begin{equation}\label{def-needlets1}
\ph_\xi(x) := \lambda_\xi^{1/2}\Phi_j(x, \xi)
\quad\mbox{and}\quad
\psi_\xi(x) := \lambda_\xi^{1/2}\Psi_j(x, \xi),
\qquad \xi \in \cX_j.
\end{equation}
Notice that for $\xi\in\cX_1$, we have
$\ph_\xi(x)=\ha(1)\PP_1(x, \xi)$ and $\psi_\xi(x)=\hb(1)\PP_1(x, \xi)$,
but $\PP_1(\cdot, \xi)\equiv 0$ if and only if $\xi=0$.
So, to prevent $\psi_\xi\equiv 0$ and $\psi_\xi\equiv 0$ for $\xi\in\cX_1$,
we (may) assume that $0\notin \cX_1$.

We set $\cX := \cup_{j = 0}^\infty \cX_j$,
where equal points from different levels $\cX_j$ are considered
as distinct elements of $\cX$, so that $\cX$ can be used as an index set.
We define the {\em analysis} and {\em synthesis} needlet systems
$\Phi$ and $\Psi$ by
\begin{equation}\label{def-needlets2}
\Phi:=\{\ph_\xi\}_{\xi\in\cX}, \quad \Psi:=\{\psi_\xi\}_{\xi\in\cX}.
\end{equation}

Estimate (\ref{eq:main-estI}) yields the rapid decay of needlets,
namely, for $x\in B^d$,
\begin{equation}\label{local-Needlets2}
|\Phi_j(\xi, x)|, |\Psi_j(\xi, x)| \le \frac{c_k 2^{jd}}
{\sqrt{\WW(2^j; \xi)}\sqrt{\WW(2^j; x)}(1+2^{j}d(\xi, x))^k} \quad
\forall k,
\end{equation}
and hence
\begin{equation}\label{local-Needlets21}
|\ph_\xi(x)|, |\psi_\xi(x)|
\le \frac{c_k 2^{jd/2}}{\sqrt{\WW(2^j;x)}(1+2^{jd}d(\xi, x))^k}
\quad \forall k.
\end{equation}
Note that on account of \eqref{eq:useful} $x$ in the term $\sqrt{\WW(2^j; x)}$
in \eqref{local-Needlets21} can be replaced by~$\xi$.

The needlets are {\rm Lip 1} functions in the following sense:
Let $\xi\in\cX_j$, $j\ge 0$, $c^*>0$, and $\omega\in B^d$. Then for each $x\in B_\omega(c^*2^{-j})$
\begin{equation}\label{Lip-Needlets}
|\ph_\xi(x)-\ph_\xi(\omega)|, |\psi_\xi(x)-\psi_\xi(\omega)|
\le \frac{c_k 2^{j(d/2+1)}d(\omega, x)}{\sqrt{\WW(2^j;\xi)}(1+2^{jd}d(\xi, \omega))^k}
\quad \forall k.
\end{equation}
This estimate follows readily from \eqref{Lip}.

We shall need estimates of the norms of the needlets.
By \eqref{eq:LpUB}, \eqref{eq:est-Lp-norm}, and since $0\notin \cX_1$,
we have for $0<p\le\infty$,
\begin{equation}\label{norm-Needlets}
\|\ph_\xi\|_\Lp \sim \|\psi_\xi\|_\Lp \sim \|\tONE_{R_\xi}\|_\Lp
\sim \Big(\frac{2^{jd}}{\WW(2^j; \xi)}\Big)^{1/2-1/p},
\quad \xi\in\cX_j.
\end{equation}
Furthermore, there exist constants $c^*, c >0$ such that
\begin{equation}\label{norm-Needlets2}
\|\ph_\xi\|_{L_\infty(B_\xi(c^*2^{-j}))},\;
\|\psi_\xi\|_{L_\infty(B_\xi(c^*2^{-j}))}
\ge c \Big(\frac{2^{jd}}{\WW(2^j; \xi)}\Big)^{1/2},
\quad \xi\in\cX_j.
\end{equation}
The proof of (\ref{norm-Needlets2}) is given in \S\ref{proofsA}.
Notice that if $\ha$, $\hb$ are real valued, then Lemma~\ref{L2-lower-bound} bellow
yields
$$
|\ph_\xi(\xi)|, |\psi_\xi(\xi)|\ge c \Big(\frac{2^{jd}}{\WW(2^j; \xi)}\Big)^{1/2},
\quad \xi\in\cX_j.
$$


Our first step in implementing needlets is to establish needlet decompositions
of $\cD'$ and $\Lpp$.


\begin{proposition}\label{prop:needlet-rep}
$(a)$ For any $f \in \cD'$,
\begin{equation}\label{Needle-rep}
f = \sum_{j=0}^\infty
\Psi_j*\overline{\Phi}_j*f
\quad\mbox{in} ~ \cD'
\end{equation}
and
\begin{equation}\label{needlet-rep1}
f = \sum_{\xi \in \cX}
\langle f, \ph_\xi\rangle \psi_\xi
\quad\mbox{in} ~ \cD'.
\end{equation}

$(b)$ For $f \in \Lpp$, $1\le p \le \infty$,
$(\ref{Needle-rep})-(\ref{needlet-rep1})$ hold in $\Lpp$.
Moreover, if $1 < p < \infty$, then the convergence in
$(\ref{Needle-rep})-(\ref{needlet-rep1})$ is unconditional.
\end{proposition}


\begin{proof}
By Definition~\ref{Def-conv} and
\eqref{def.Phi-j} we have, for $f\in\cD'$,
\begin{equation}\label{Phi-f}
\overline{\Phi}*f = \sum_{\nu=0}^{2^j}\overline{\ha\Big(\frac{\nu}{2^{j-1}}\Big)}\PP_\nu*f
\end{equation}
and using Lemma~\ref{lem:prop-conv} and that $\PP_\nu*\PP_\nu(\cdot, y)=\PP_\nu(\cdot, y)$
\begin{equation}\label{Psi-Phi-f}
\Psi*\overline{\Phi}*f
= \sum_{\nu=0}^{2^j}\overline{\ha\Big(\frac{\nu}{2^{j-1}}\Big)}\hb\Big(\frac{\nu}{2^{j-1}}\Big)\PP_\nu*f.
\end{equation}
Then \eqref{Needle-rep} follows from the above, \eqref{part-unity}, and Lemma~\ref{lem:dec-D1}.

Note that $\Psi_j(x, y)\overline{\Phi(y, z)}$ belongs to $\Pi_{2^{j+1}-1}$
as a function of $y$ and, therefore, employing the cubature formula from
Proposition~\ref{prop:cubature} we get
\begin{align*}
\Psi_j*\overline{\Phi_j(\cdot, z)}
&=\int_{B^d}\Psi_j(x, y)\overline{\Phi(y, z)}\W(y)dy\\
&=\sum_{\xi\in\cX_j}\lambda_\xi \Psi_j(x, \xi)\overline{\Phi(\xi, z)}
=\sum_{\xi\in\cX_j}\psi_\xi(x)\overline{\ph_\xi(z)},
\end{align*}
which leads to
$$
\Psi_j*\overline{\Phi}_j*f=\sum_{\xi\in\cX_j}\langle f, \ph_\xi\rangle \psi_\xi.
$$
Combining this with \eqref{Needle-rep} yields \eqref{needlet-rep1}.

The convergence of \eqref{Needle-rep} and \eqref{needlet-rep1} in
$\Lpp$ for $f\in \Lpp$ follows in a similar fashion (see also
\cite[Proposion 3.1]{KPX}). The unconditional convergence in
$\Lpp$, $1<p<\infty$, follows by
Theorem~\ref{thm:Fnorm-equivalence} and
Proposition~\ref{prop:ident} below.
\end{proof}

\section{Weighted Triebel-Lizorkin  spaces on \boldmath $B^d$}\label{Tri-Liz}
\setcounter{equation}{0}

Following the general idea of using spectral decompositions (see e.g. \cite{Pee, T1}),
we next employ orthogonal polynomials to introduce
weighted Triebel-Lizorkin spaces on $B^d$.
To this end we define a sequence of kernels $\{\Phi_j\}$ by
\begin{equation}\label{def-Phi-j}
\Phi_0(x, y) := 1
\quad\mbox{and}\quad
\Phi_j(x, y) := \sum_{\nu=0}^\infty \ha
\Big(\frac{\nu}{2^{j-1}}\Big)\PP_\nu(x,y), ~~~j\ge 1,
\end{equation}
where $\{\PP_\nu(x,y)\}$ are from \eqref{def-Pn}-\eqref{compactPn} and
$\ha$ obeys the conditions
\begin{align}
&\quad \ha\in C^\infty[0, \infty),
\quad \supp \, \ha \subset [1/2, 2], \label{ha1}\\
&\quad  |\ha(t)|>c>0, \quad \text{if } t \in [3/5, 5/3].\label{ha2}
\end{align}


\begin{definition}
Let  $s, \r \in \R$, $0<p<\infty$, and $0<q\le\infty$. Then the
weighted Triebel-Lizorkin space $\Fsrpq:=\Fsrpq(\W)$ is defined as
the set of all $f\in \cD'$ such that
\begin{equation}\label{Tri-Liz-norm}
\|f\|_{\Fsrpq}:=\Big\|\Big(\sum_{j=0}^{\infty}
\Big[2^{sj}\WW(2^j;\cdot)^{-\r/d}|\Phi_j*f(\cdot)|\Big]^q\Big)^{1/q}\Big\|_{\Lp} <\infty
\end{equation}
with the usual modification when $q=\infty$.
\end{definition}


Observe that the above definition is independent of the choice
of $\ha$ as long as it satisfies
(\ref{ha1})-(\ref{ha2}) (see Theorem~\ref{thm:Fnorm-equivalence} below).


\begin{proposition}\label{Fsrpq-embedding}
For all $s,\r \in \R$, $0<p<\infty$, and $0<q\le\infty$,
$\Fsrpq$ is a quasi-Banach space which is continuously embedded in $\cD'$.
\end{proposition}

\begin{proof}The completeness of the space $\Fsrpq$
follows easily (see e.g. \cite{T1}, p. 49) by the continuous
embedding of $\Fsrpq$ in $\cD'$, which we establish next.

Let $\{\Phi_j\}$ be the kernels from the definition of $\Fsrpq$ with $\ha$ obeying
(\ref{ha1})-(\ref{ha2}) that are the same as (\ref{ha-hb1})-(\ref{ha-hb2}).
As already indicated there exists a function $\hb$ satisfying (\ref{ha-hb1})-(\ref{ha-hb3}).
We use this function to define $\{\Psi_j\}$ as in (\ref{def-Psi-j}).
Then by Proposition~\ref{prop:needlet-rep}
$f=\sum_{j=0}^\infty \Psi_j*\Phi_j*f$ in $\cD'$ and hence
$$
\langle f, \phi\rangle=\sum_{j=0}^\infty \langle \Psi_j*\overline{\Phi}_j*f, \phi\rangle,
\quad \phi\in\cD.
$$
We now employ \eqref{Phi-f}-\eqref{Psi-Phi-f} and the Cauchy-Schwarz inequality
to obtain, for $j\ge 2$,
\begin{align*}
&|\langle \Psi_j\ast\overline{\Phi}_j\ast f, \phi\rangle|^2
=\Big|\sum_{\nu=2^{j-2}+1}^{2^j}
\overline{\ha\Big(\frac{\nu}{2^{j-1}}\Big)}\hb\Big(\frac{\nu}{2^{j-1}}\Big)
\langle \proj_\nu f, \proj_\nu \phi\rangle\Big|^2\\
&\qquad \le\sum_{\nu=2^{j-2}+1}^{2^j}
\Big|\ha\Big(\frac{\nu}{2^{j-1}}\Big)\Big|^2\|\proj_\nu f\|_2^2
\sum_{\nu=2^{j-2}+1}^{2^j}
\Big|\hb\Big(\frac{\nu}{2^{j-1}}\Big)\Big|^2\|\proj_\nu \phi\|_2^2\\
&\qquad\le 2^j\|\Phi_j*f\|_2^2
\max_{2^{j-2}< \nu\le 2^j}\|\proj_\nu \phi\|_2^2.
\end{align*}
Using inequality \eqref{norm-relation} we get
\begin{align*}
\|\Phi_j*f\|_2 \le c2^{j(d+2\mu)/p}\|\Phi_j*f\|_p
\le c2^{j((d+2\mu)/p+2\mu|\r|/d-s)}\|2^{sj}\WW(2^j;\cdot)^{-\r/d}\Phi_j*f(\cdot)\|_p.
\end{align*}
From the above estimates we infer
\begin{align*}
|\langle \Psi_j\ast\overline{\Phi}_j\ast f, \phi\rangle|
\le c 2^{-j}\|f\|_{\Fsrpq} 2^{jk} \max_{2^{j-2}< \nu\le 2^j}\|\proj_\nu f\|_2
\le c 2^{-j}\|f\|_{\Fsrpq}\cN_k(\phi)
\end{align*}
for $k\ge (d+2\mu)/p+2\mu|\r|/d +3/2-s$.
A similar estimate trivially holds for $j=0, 1$.
Summing up we get
$$
|\langle f, \phi\rangle| \le c\|f\|_{\Fsrpq}\cN_k(\phi),
$$
which completes the proof.
\end{proof}

\medskip

As a companion to $\Fsrpq$ we now introduce the sequence spaces $\fsrpq$.
Here we assume that $\{\cX_j\}_{j=0}^\infty$ is a sequence of almost
uniformly $\eps_j$-distributed points on $B^d$
$(\eps_j:=c^\diamond 2^{-j})$
with associated neighborhoods $\{R_\xi\}_{\xi\in\cX_j}$,
given by Proposition~\ref{prop:cubature}.
Just as in the definition of needlets in \S\ref{def-needlets},
we set $\cX:=\cup_{j\ge 0} \cX_j$.


\begin{definition}
Suppose $s, \r \in \R$, $0<p<\infty$, and $0<q\le\infty$. Then~$\fsrpq$
is defined as the space of all complex-valued sequences
$h:=\{h_{\xi}\}_{\xi\in \cX}$ such that
\begin{equation}\label{def-f-space}
\norm{h}_{\fsrpq} :=\nnorm{\biggl(\sum_{j=0}^\infty
2^{sjq}\sum_{\xi \in \cX_j}
(|h_{\xi}|\WW(2^j;\xi)^{-\r/d}\tONE_{R_\xi}(\cdot))^q\biggr)^{1/q}}_{\Lp} <\infty
\end{equation}
with the usual modification for $q=\infty$. Recall that
$\tONE_{R_\xi}:=m(R_\xi)^{-1/2}\ONE_{R_\xi}$.
\end{definition}

In analogy to the classical case on $\R^d$ we introduce ``analysis" and ``synthesis" operators by
\begin{equation}\label{anal_synth_oprts}
S_\varphi: f\rightarrow \{\ip{f, \varphi_\xi}\}_{\xi \in \cX}
\quad\text{and}\quad
T_\psi: \{h_\xi\}_{\xi \in \cX}\rightarrow \sum_{\xi\in \X}h_\xi\psi_\xi.
\end{equation}


We now give our main result on weighted Triebel-Lizorkin spaces.


\begin{theorem}\label{thm:Fnorm-equivalence}
Let $s, \r\in \R$, $0< p< \infty$ and $0<q\le \infty$. Then the operators
$S_\varphi:\Fsrpq\rightarrow\fsrpq$ and $T_\psi:\fsrpq\rightarrow \Fsrpq$
are bounded and $T_\psi\circ S_\varphi=Id$ on $\Fsrpq$.
Consequently, $f\in \Fsrpq$ if and only if
$\{\ip{f, \varphi_\xi}\}_{\xi \in \cX}\in \fsrpq$.
Furthermore,
\begin{align}\label{Fnorm-equivalence-1}
\norm{f}_{\Fsrpq} &
\sim  \norm{\{\ip{f,\varphi_\xi}\}}_{\fsrpq}
\sim \nnorm{\biggl(\sum_{j=0}^\infty 2^{sjq}
\sum_{\xi\in \cX_j}(|\ip{f, \varphi_\xi}|\WW(2^j; \xi)^{-\r/d}|\psi_\xi(\cdot)|)^q\biggr)^{1/q}}_{\Lp}.
\end{align}
In addition, the definition of $\Fsrpq$ is independent of the particular selection of
$\ha$ satisfying $(\ref{ha1})$--$(\ref{ha2})$.
\end{theorem}

The proof of this theorem relies on several lemmas whose proofs are given
in Section \ref{proofsB}.
In the following we assume that $\{\Phi_j\}$ are from the definition
of weighted Triebel-Lizorkin spaces,
while $\{\ph_\xi\}_{\xi\in\cX}$ and $\{\psi_\xi\}_{\xi\in\cX}$
are needlet systems defined as in (\ref{def-needlets1})
with no connection between the functions $\ha$'s from (\ref{def-Phi-j}) and (\ref{def.Phi-j}).


\begin{lemma}\label{lem:Phi*psi}
For any $k > 0$ there exists a constant $c_k>0$ such that
\begin{equation}\label{Phi*psi}
|\Phi_j\ast \psi_\xi(x)| \le c_k\frac{2^{jd/2}}{\sqrt{ \WW(2^j;x)
}(1+2^jd(x, \xi))^k}, \quad \xi\in\cX_\nu, \quad  j-1\le \nu\le j+1,
\end{equation}
and $\Phi_j\ast \psi_\xi\equiv 0$ for $\xi\in\cX_\nu$ if $\nu\ge j+2$ or
$\nu\le j-2$. Here $\cX_\nu:=\emptyset$ if $\nu < 0$.
\end{lemma}


\begin{lemma}\label{lem:Max-needl}
For any $t>0$ and $\xi\in \cX_j,$ $j\ge 0$,
\begin{equation}\label{Max-needl1}
|\ph_\xi(x)|, |\psi_\xi(x)| \le c (\cM_t\tONE_{R_\xi})(x), \quad
x\in B^d, \quad\mbox{and}
\end{equation}
\begin{equation}\label{Max-needl2}
\tONE_{R_\xi}(x) \le c(\cM_t \ph_\xi)(x), c(\cM_t \psi_\xi)(x),
\quad x\in B^d.
\end{equation}
\end{lemma}


\begin{definition}\label{def-h-star}
For any set of complex numbers
$\{h_{\xi}\}_{\xi\in \cX_j}$ $(j\ge 0)$ we define
\begin{equation}\label{def.h-star}
h^{\ast}_\xi:=\sum_{\eta\in
\cX_j}\frac{|h_\eta|}{(1+2^jd(\eta,\xi))^\sigma},
\quad \xi\in \cX_j,
\end{equation}
where $\sigma>1$ is a sufficiently large constant that will be selected later on.
\end{definition}


\begin{lemma}\label{lem:weak_inequality}
Let $P\in \Pi_{2^j}$, $j\ge 0$, and denote
$a_\xi:=\max_{x\in R_\xi}|P(x)|$ for $\xi\in\cX_j$.
There exists $r\ge 1$, depending only on $\sigma$, $\mu$, and $d$
such that if
$$
b_\xi:=\max\{\min_{x\in R_\eta }|P(x)|:\eta\in \cX_{j+r},
R_\xi\cap R_\eta\ne\emptyset \}, \quad \xi\in\cX_j,
$$
then
\begin{equation}\label{weak_inequality}
a_\xi^\ast\sim b_\xi^\ast
\end{equation}
with constants of equivalence independent of $P$, $j$, and $\xi$.
\end{lemma}


\begin{lemma}\label{lem:disc-max}
Assume $t>0$, $\gamma \in \R$, and let $\{b_\xi\}_{\xi\in \cX_j}$
$(j\ge 0)$ be a set of complex numbers. Also, let $\sigma$ in the
definition $(\ref{def.h-star})$ of $b_\xi^*$ obey
$\sigma>d+(d+2\mu)/t+2\mu|\gamma|$.
Then for any $\xi \in \cX_j$
\begin{equation}\label{disc-max}
b_\xi^{\ast}\WW(2^j; \xi)^{\gamma}\ONE_{R_\xi}(x)
\le c \cM_t\Big(\sum_{\eta\in \cX_j}|b_\eta|\WW(2^j; \eta)^{\gamma}\ONE_{R_\eta}(\cdot)\Big)(x),
\quad x\in R_\xi.
\end{equation}
\end{lemma}


\noindent {\em Proof of Theorem $\ref{thm:Fnorm-equivalence}$.}
Choose $0<t<\min\{p,q\}$ and let $\sigma$ in
Definition~\ref{def-h-star} obey $\sigma
>d+(d+2\mu)/t+2\mu|\r|/d$. Now, choose $k\ge\sigma+2\mu|\r|/d$.
Observe first that the right-hand side equivalence in
(\ref{Fnorm-equivalence-1}) follows immediately from Lemma
\ref{lem:Max-needl} and the maximal inequality (\ref{max-ineq}).

Let $\{\Phi_j\}$  be a sequences of kernels as in
the definition of weighted Triebel-Lizorkin spaces,
i.e. $\Phi_j$ is defined by (\ref{def-Phi-j})
with $\ha$ satisfying (\ref{ha1})-(\ref{ha2}),
the same as (\ref{ha-hb1})-(\ref{ha-hb2}).
As already mentioned, there exists a function $\hb$ satisfying
(\ref{ha-hb1})-(\ref{ha-hb2}) such that (\ref{ha-hb3}) holds.
Let $\Psi_j$ be defined by (\ref{def-Psi-j}) with this $\hb$.
In addition, let $\{\ph_\xi\}_{\xi\in\cX}$ and $\{\psi_\xi\}_{\xi\in\cX}$
be the associated needlet systems defined as in (\ref{def-needlets1})
using these $\ha$ and $\hb$.

Exactly in the same way, let $\{\wt\Phi_j\}$ and $\{\wt\Psi_j\}$ be two sequences
of kernels defined as above using completely different functions $\ha$ and $\hb$.
Also, assume that $\{\wt\ph_\xi\}$, $\{\wt\psi_\xi\}$
are the associated needlet systems, defined as in
(\ref{def.Phi-j})-(\ref{def-needlets1}).
As a result, we have two completely different systems of kernels and associated needlet systems.

Let us first prove the boundedness of the operator
$T_{\wt\psi}:\fsrpq\rightarrow\Fsrpq$,
defined similarly as in~(\ref{anal_synth_oprts}) with $\{\psi_\xi\}$
replaced by $\{\widetilde{\psi}_\xi\}$.
Here we assume that space $\Fsrpq$ is defined by $\{\Phi_j\}$.
Let $h:=\{h_\xi\}_{\xi\in \cX}$ be an arbitrary finitely supported sequence and
$f:=\sum_{\xi}h_\xi\widetilde{\psi}_\xi$.
Using Lemma~\ref{lem:Phi*psi} we have, for $x\in B^d$,
\begin{align*}
|\Phi_j\ast f(x)| &=\Big|\sum_{\xi\in \cX}  h_\xi
\Phi_j\ast\wt\psi_{\xi}(x)\Big| \le \sum_{j-1\le \nu\le
j+1}\sum_{\xi\in \cX_\nu} |h_\xi |
|\Phi_j\ast\wt\psi_{\xi}(x)|\\
&\le c2^{jd/2} \sum_{j-1\le \nu\le j+1}\sum_{\xi\in \cX_\nu} \frac{
|h_\xi| }{\sqrt{ \WW(2^\nu;x)}(1+2^\nu d(\xi, x))^k}.
\end{align*}
For $\eta\in\cX_j$, denote
$\Gamma_\eta:=\{w\in\cX_{j-1}\cup\cX_j\cup\cX_{j+1}: R_w\cap R_\eta\ne \emptyset\}$.
Here $\cX_{-1}:=\emptyset$.
Note first that $\# \Gamma_\eta \le c$.
Secondly, for $x\in R_\eta$ and $w\in\Gamma_\eta$, we have $d(x, w) \le c2^{-j}$ and
using inequality \eqref{eq:useful}
$$
\WW(2^j; x)^{-\r/d} \le c\WW(2^j; w)^{-\r/d}
\le c\WW(2^j; \xi)^{-\r/d}(1+2^jd(\xi, \omega))^{2\mu|\r|/d}.
$$
We use the above estimates to obtain, for $x\in R_\eta$,
\begin{align*}
&\WW(2^j; x)^{-\r/d}|\Phi_j\ast f(x)|\\
&\qquad\qquad \le c2^{jd/2} \sum_{j-1\le \nu\le j+1}
\sum_{\omega\in\Gamma_\eta\cap\cX_\nu}\sum_{\xi\in \cX_\nu}
\frac{|h_\xi|\WW(2^j; \xi)^{-\r/d}\ONE_{R_\omega}(x)}
{\sqrt{\WW(2^\nu;\omega)}(1+2^\nu d(\xi, \omega))^{k-2\mu|\r|/d}}\\
&\qquad\qquad \le c2^{jd/2}
\sum_{\omega\in\Gamma_\eta}
\frac{H_\omega^*\ONE_{R_\omega}(x)}{\sqrt{\WW(2^j;\omega)}}
\le c\sum_{\omega\in\Gamma_\eta}H_\omega^*\tONE_{R_\omega}(x),
\end{align*}
where
$H_\omega:= h_\omega\WW(2^j; \omega)^{-\r/d}$.
Here we used that $k-2\mu|\r|/d \ge \sigma$
and \eqref{eq:mR}.
We insert the above in (\ref{Tri-Liz-norm}) and use Lemma~\ref{lem:disc-max} (with $\gamma=0$)
and the maximal inequality (\ref{max-ineq}) to obtain
\begin{equation}
\begin{aligned}\label{F-less-f}
\norm{f}_{\Fsrpq}
&\le c\nnorm{\Bigr(\sum_{j=0}^\infty \Big[2^{sj}\sum_{\eta\in \cX_j}\sum_{\omega\in \Gamma_\eta}
H_\omega^*\tONE_{R_\omega}(\cdot)\Big]^q\Big)^{1/q}}_{\Lp}\\
&\le c\nnorm{\Bigr(\sum_{j=0}^\infty \Big[2^{sj}\sum_{\xi\in \cX_j}
H_\xi^*\tONE_{R_\xi}(\cdot)\Big]^q\Big)^{1/q}}_{\Lp}\\
&\le c\nnorm{\Bigr(\sum_{j=0}^\infty \Big[\cM_t\Big(\sum_{\xi\in \cX_j}
2^{sj}|H_\xi|\tONE_{R_\xi}\Big)(\cdot)\Big]^q\Big)^{1/q}}_{\Lp}\\
&\le c\nnorm{\Bigr(\sum_{j=0}^\infty \Big[\sum_{\xi\in \cX_j}
2^{sj}|H_\xi|\tONE_{R_\xi}(\cdot)\Big]^q\Big)^{1/q}}_{\Lp}
\le c\norm{\{h_\xi\} }_{\fsrpq},
\end{aligned}
\end{equation}
where in  the second inequality  above we used that
$\#\Gamma_\eta\le c$.
This establishes the desired result for finitely supported
sequences. Using  the continuous  embedding of $\Fsrpq$ in $\cD'$
(Proposition~\ref{Fsrpq-embedding}) and the density of finitely
supported sequences in $\fsrpq$ it follows from  (\ref{F-less-f})
that for every $h\in \fsrpq$, $T_{\wt\psi} h:=\sum_{\xi\in
\cX}h_\xi\wt\psi_\xi$ is a well defined distribution in $\cD'$.
Then a standard density argument shows that
$T_{\wt\psi}:\fsrpq\rightarrow\Fsrpq$ is bounded.


Assume now that the space $\Fsrpq$ is defined in terms of $\{\overline{\Phi}_j\}$
in place of $\{\Phi_j\}$.
Using this definition we shall prove the boundedness of the operator
$S_\varphi:\Fsrpq\rightarrow \fsrpq$.

Let $f\in \Fsrpq$. Then $\overline{\Phi}_j\ast f\in \Pi_{2^j}$.
For $\xi \in \cX_j$, we define
$$
a_\xi:=\max_{x\in R_\xi}|\overline{\Phi}_j\ast f(x)|, \quad
b_\xi:=\max\{\min_{x\in R_\eta }|\overline{\Phi}_j\ast f(x)|:\eta\in
\cX_{j+r}, R_\xi\cap R_\eta\ne\emptyset \},
$$
where $r\ge 1$ is from Lemma \ref{lem:weak_inequality}.
Then by the same lemma $a_\xi^*\sim b_\xi^*$.
Hence, using \eqref{eq:mR},
$$
|\langle f, \ph_\xi
\rangle|=\lambda_\xi^{1/2}|\overline{\Phi}_j*f(\xi)| \le c
m(R_\xi)^{1/2}a_\xi \le c m(R_\xi)^{1/2}a_\xi^* \le c
m(R_\xi)^{1/2}b_\xi^*.
$$
From this, recalling that
$\tONE_{R_\xi}:=m(R_\xi)^{-1/2}\ONE_{R_\xi}$, we get
\begin{equation}\label{asdfasf}
\begin{aligned}
\norm{\{\ip{f, \varphi_\xi}\}}_{\fsrpq}
&=\nnorm{\Bigl(\sum_{j=0}^\infty 2^{jsq}
\sum_{\xi\in \cX_j}
[|\ip{f, \varphi_\xi}|\WW(2^j; \xi)^{-\r/d}\tONE_{R_\xi}(\cdot)]^q\Bigr)^{1/q}}_{\Lp}\\
&\le c\nnorm{\Bigl(\sum_{j=0}^\infty 2^{jsq}
\sum_{\xi\in \cX_j}[b_\xi^*\WW(2^j; \xi)^{-\r/d}
\ONE_{R_\xi}(\cdot)]^q\Bigr)^{1/q}}_{\Lp}\\
&\le c\nnorm{\Bigl(\sum_{j=0}^\infty 2^{jsq}
\Big[\cM_t\Big(\sum_{\xi\in \cX_j}
b_\xi\WW(2^j; \xi)^{-\r/d}\ONE_{R_\xi}(\cdot)\Big)(\cdot)\Big]^q\Bigr)^{1/q}}_{\Lp}\\
&\le c\nnorm{\Bigl(\sum_{j=0}^\infty2^{jsq}
\Big[\sum_{\xi\in \cX_j}
b_\xi \WW(2^j; \xi)^{-\r/d}\ONE_{R_\xi}(\cdot)\Big]^q\Bigr)^{1/q}}_{\Lp}.
\end{aligned}
\end{equation}
Here for the second inequality above we used Lemma~\ref{lem:disc-max}
and for the third one the maximal inequality (\ref{max-ineq}).

Denote $m_\eta:=\min_{x\in R_\eta}|\overline{\Phi}_j\ast f(x)|$ for
$\eta\in \cX_{j+r}$
and
$$
\cX_{j+r}(\xi):=\{w\in \cX_{j+r}:R_w\cap R_\xi \ne \emptyset\}
\quad\mbox{for $\xi\in\cX_j$}.
$$
Evidently $\#\cX_{j+r}(\xi)\le c(r,d)$.
Further, for $w, \eta\in \cX_{j+r}(\xi)$ we have
$d(w,\eta)\le c2^{-j}$ and hence
$$
m_w\le c\frac{m_w}{(1+2^{j+r}d(w,\eta))^\sigma}\le cm^\ast_\eta,
\quad c=c(r, \sigma, d).
$$
Therefore, for any $\eta\in \cX_{j+r}(\xi)$,
$
b_\xi=\max_{w\in\cX_{j+r}(\xi)}m_w \le cm^\ast_\eta.
$
and hence
\begin{equation}\label{b-xi-m-eta}
b_\xi\ONE_{R_\xi}\le \sum_{\eta\in \cX_{j+r}(\xi)}  m^\ast_\eta \ONE_{R_\eta}.
\end{equation}
Clearly,
$\WW(2^j; \xi) \sim \WW(2^{j+r}; \eta)$ for $\eta\in\cX_{j+r}(\xi)$.
This along with \eqref{b-xi-m-eta} leads to
\begin{equation}\label{W-b-xi-m-eta}
b_\xi\WW(2^j;\xi)^{-\r/d}\ONE_{R_\xi}
\le c\sum_{\eta\in \cX_{j+r}(\xi)}  m^*_\eta\WW(2^{j+r}; \eta)^{-\r/d}\ONE_{R_\eta}.
\end{equation}
Using  this estimate in (\ref{asdfasf}) we get
\begin{equation*}
\begin{aligned}
\norm{\{\ip{f, \varphi_\xi}\}}_{\fsrpq}
&\le c\nnorm{\Bigl(\sum_{j=0}^\infty 2^{jsq}
\Bigl(\sum_{\eta\in \cX_{j+r}}  m^*_\eta\WW(2^{j+r}; \eta)^{-\r/d}
\ONE_{R_\eta}(\cdot) \Bigr)^q\Bigr)^{1/q}}_{\Lp}\\
&\le c\nnorm{\Bigl(\sum_{j=0}^\infty 2^{jsq}
\Big[\cM_t\Big(\sum_{\eta\in \cX_{j+r}}m_\eta\WW(2^{j+r}; \eta)^{-\r/d}
\ONE_{R_\eta}\Big)(\cdot)\Big]^q\Bigr)^{1/q}}_{\Lp}\\
&\le c\nnorm{\Bigl(\sum_{j=0}^\infty \Bigl(2^{js}
\sum_{\eta\in \cX_{j+r}} m_\eta\WW(2^{j+r}; \eta)^{-\r/d}\ONE_{R_\eta}(\cdot) \Bigr)^q\Bigr)^{1/q}}_{\Lp}\\
&\le c\nnorm{\Bigl(\sum_{j=0}^\infty (2^{js}
 \WW(2^{j}; \cdot)^{-\r/d} |\overline{\Phi}_j*f(\cdot)|)^q\Bigr)^{1/q}}_{\Lp}
 =c\|f\|_{\Fsrpq}.
\end{aligned}
\end{equation*}
Here for first inequality we used that $\# \cX_{j+r}(\xi) \le c$,
for the second inequality we used Lemma~\ref{lem:disc-max},
and for third one the maximal inequality \eqref{max-ineq}.
We also use that $\WW(2^{j+r}; \eta)\sim \WW(2^{j}; x)$ if $x\in \R_\eta$, $\eta\in\cX_{j+r}$.
Thus the boundedness of $S_\varphi:\Fsrpq\rightarrow \fsrpq$ is established.

The identity $T_\psi\circ S_\varphi=Id$ follows by Proposition~\ref{prop:needlet-rep}.

It remains to show that $\Fsrpq$ is independent of the particular selection
of $\ha$ in the definition of  $\{\Phi_j\}$.
Denote by $\|\cdot\|_{\Fsrpq(\Phi)}$ the F-norm defined by $\{\Phi_j\}$.
Then by the above proof it follows that
$$
\|f\|_{\Fsrpq(\Phi)}
\le c\|\{\langle f, \wt\ph_\xi\rangle\}\|_{\fsrpq}
\quad\mbox{and}\quad
\|\{\langle f, \ph_\xi\rangle\}\|_{\fsrpq}
\le c\|f\|_{\Fsrpq(\overline{\Phi})}
$$
and hence
$$
\|f\|_{\Fsrpq(\Phi)}
\le c\|\{\langle f, \wt\ph_\xi\rangle\}\|_{\fsrpq}
\le c\|f\|_{\Fsrpq(\overline{\wt\Phi})}.
$$
Now the desired independence follows by interchanging the roles of
$\{ \Phi_j\}$,$\{\wt\Phi_j\}$, and  their complex conjugates.
$\qed$

\medskip

In a sense the spaces $\F sspq$ are more natural than the spaces
$\Fsrpq$ with $\r\ne s$ since they scale (are embedded)
``correctly" with respect to the smoothness index $s$.

%
\begin{proposition}\label{F-embedding}
Let $0<p<p_1<\infty$, $0<q, q_1\le\infty$, and
$-\infty<s_1<s<\infty$. Then we have the continuous embedding
\begin{equation}\label{F-embed}
\F sspq \subset \F {s_1}{s_1}{p_1}{q_1} \quad\mbox{if}\quad
s/d-1/p=s_1/d-1/p_1.
\end{equation}
\end{proposition}

The proof of this embedding result can be carried out similarly as
in the classical case on $\RR^n$ using inequality
(\ref{norm-relation2}) and Theorem~\ref{thm:Fnorm-equivalence}
(see e.g. \cite{T1}, page 129). It~will be omitted.

\medskip


Finally, we would like to link the weighted Triebel-Lizorkin spaces $\Fsrpq$
to $\Lpp$ and weighted potential space
(generalized weighted Sobolev space) on $B^d$.

We define the weighted potential space
$H_p^s := H_p^s(\W)$,
$s >0$, $1\le p \le \infty$,  on $B^d$
as the set of all $f\in \cD'$ such that
\begin{equation}\label{def-H}
\|f\|_{H_p^s}:=
\Big\|\sum_{n=0}^\infty (n+1)^s \proj_n f\Big\|_\Lp
< \infty,
\end{equation}
where
$\proj_n f:=\PP_n*f$.

We have the following identification of certain weighted Triebel-Lizorkin spaces.


\begin{proposition}\label{prop:ident}
We have
$$ 
\F s0p2 \sim H_p^s,
\quad s >0, ~1 < p < \infty,
$$ 
and
$$ 
\F 00p2\sim \Lpp,
\quad 1 < p < \infty,
$$ 
with equivalent norms.
Consequently, for any $f\in\Lpp$, $1<p<\infty$,
$$
\|f\|_p
\sim \nnorm{\biggl(\sum_{j=0}^\infty
\sum_{\xi\in \cX_j}(|\ip{f, \varphi_\xi}||\psi_\xi(\cdot)|)^2\biggr)^{1/2}}_{\Lp}.
$$
\end{proposition}

The proof of this proposition uses the multipliers from \cite[Theorem 5.2]{DX}
and can be carried out exactly as in the case of spherical
harmonic expansions in \cite[Proposition~4.3]{NPW2}.
We omit it.

\section{Weighted Besov spaces on \boldmath $B^d$}\label{Besov}
\setcounter{equation}{0}

For the definition of weighted Besov spaces on $B^d$ we use the sequence
of kernels $\{\Phi_j\}$ defined in (\ref{def-Phi-j}) with $\ha$
obeying (\ref{ha1})-(\ref{ha2})
(see \cite{Pee, T1} for the general idea of using spectral decompositions).


\begin{definition}
Let $s, \r\in \R$ and $0<p,q \le \infty$. The weighted Besov space
$\Bsrpq := \Bsrpq(\W)$ is defined
as the set of all $f \in \cD'$ such that
\begin{equation}\label{def-Besov-sp}
\|f\|_{\Bsrpq} :=
\Big(\sum_{j=0}^\infty \Big(2^{s j}
\|\WW(2^j; \cdot)^{-\r/d}\Phi_j*f(\cdot)\|_{\Lp}\Big)^q\Big)^{1/q}
< \infty,
\end{equation}
where the $\ell_q$-norm is replaced by the sup-norm if $q=\infty$.
\end{definition}

Observe that as in the case of weighted Triebel-Lizorkin spaces the above definition
is independent of the particular choice of $\ha$ obeying  (\ref{ha1})-(\ref{ha2})
(see Theorem~\ref{thm:Bnorm-eq}).
Also, as for $\Fsrpq$ the Besov space $\Bsrpq$ is a quasi-Banach
space which is continuously embedded in $\cD'$. We skip the details.

\smallskip

We next introduce the sequence spaces $\bsrpq$
associated to the weighted Besov spaces $\Bsrpq$.
To this end, we assume that $\{\cX_j\}_{j=0}^\infty$ is a sequence of almost
uniformly $\eps_j$-distributed points on $B^d$
$(\eps_j:=c^\diamond 2^{-j})$
with associated neighborhoods $\{R_\xi\}_{\xi\in\cX_j}$,
given by Proposition~\ref{prop:cubature}.
As before we set $\cX:=\cup_{j\ge 0} \cX_j$.


\begin{definition}
Let $s, \r\in \R$ and $0<p,q \le \infty$. Then $\bsrpq$
is defined to be the space of all complex-valued sequences
$h:=\{h_{\xi}\}_{\xi\in \cX}$ such that
\begin{equation}\label{def-berpq}
\norm{h}_{\bsrpq} :=\Bigl(\sum_{j=0}^\infty
2^{j(s-d/p+d/2)q}
\Bigl[\sum_{\xi\in \cX_j}\Big(\WW(2^j;\xi)^{-\r/d+1/p-1/2}|h_\xi|\Big)^p
\Bigr]^{q/p}\Bigr)^{1/q}
\end{equation}
is finite,
with the usual modification for $p=\infty$ or $q=\infty$.
\end{definition}

Our main result in this section is the following characterization of
weighted Besov spaces, which employs the operators $S_\ph$ and
$T_\psi$ defined in (\ref{anal_synth_oprts}).


\begin{theorem}\label{thm:Bnorm-eq}
Let $s, \r\in \R$ and  $0< p,q\le \infty$.
The operators
$S_\varphi:\Bsrpq\rightarrow\bsrpq$ and
$T_\psi:\bsrpq\rightarrow \Bsrpq$
are bounded and $T_\psi\circ S_\varphi=Id$ on $\Bsrpq$.
Consequently, for $f\in \cD'$ we have that $f\in \Bsrpq$ if and
only if $\{\ip{f, \varphi_\xi}\}_{\xi \in \cX}\in \bsrpq$.
Moreover,
\begin{align}\label{Bnorm-equivalence-1}
\norm{f}_{\Bsrpq} & \sim  \norm{\{\ip{f,\varphi_\xi}\}}_{\bsrpq}
\sim \Big(\sum_{j=0}^\infty2^{sjq} \Bigl[\sum_{\xi\in \cX_j}
\Big(\WW(2^j; \xi)^{-\r/d}\norm{\ip{f,\varphi_\xi}\psi_\xi}_{\Lp}\Big)^p\Bigr]^{q/p}\Bigr)^{1/q}.
\end{align}
In addition, the definition of $\Bsrpq$ is independent of the particular
selection of $\ha$ satisfying $(\ref{ha1})$--$(\ref{ha2})$.
\end{theorem}

For the proof of this theorem we shall utilize some of the lemmas from \S\ref{Tri-Liz}
as well as the following additional lemma
whose proof is given in Section~\ref{proofsB}.




\begin{lemma}\label{l:half_shannon}
Let $0<p\le \infty$ and $\gamma\in \R$.
Then for any $P\in \Pi_{2^j}, j\ge 0$,
\begin{equation}
 \Big(\sum_{\xi\in \cX_j}\WW(2^j; \xi)^\gamma\max_{x\in R_\xi}|P(x)|^p m(R_\xi)\Big)^{1/p}
 \le c\norm{\WW(2^j; \cdot)^\gamma P(\cdot)}_{\Lp}.
\end{equation}
\end{lemma}


\noindent {\em Proof of Theorem~$\ref{thm:Bnorm-eq}$.} We first
note that the right-hand side of (\ref{Bnorm-equivalence-1})
follows immediately from (\ref{norm-Needlets}).

Just as in the proof of Theorem~\ref{thm:Fnorm-equivalence},
we assume that $\{\Phi_j\}$ are kernels defined by (\ref{def-Phi-j}),
with $\ha$ satisfying (\ref{ha1})-(\ref{ha2}).
Next, suppose $\{\Psi_j\}$ are defined by (\ref{def-Psi-j})
with $\hb$ obeying (\ref{ha-hb1})-(\ref{ha-hb3}).
Also, let $\{\ph_\xi\}_{\xi\in\cX}$ and $\{\psi_\xi\}_{\xi\in\cX}$
be the associated needlet systems defined as in (\ref{def-needlets1}).
Further, assume that $\{\wt\Phi_j\}$, $\{\wt\Psi_j\}$,
$\{\wt\ph_\xi\}$, $\{\wt\psi_\xi\}$ is a second (completely different) set of kernels and
needlets.

Our first step is to prove the boundedness of the operator
$T_{\wt\psi}:\bsrpq\rightarrow\Bsrpq$ defined as in \eqref{anal_synth_oprts}
with $\{\psi_\xi\}$ replaced by $\{\wt\psi_\xi\}$;
we assume that $\Bsrpq$ is defined by $\{\Phi_j\}$.

Pick $0<t<\min\{p,1\}$ and $k \ge 2\mu|\r|/d+\mu+(2\mu+d)/t$.
Let $h=\{h_\xi\}_{\xi\in\cX_j}$ be a finitely supported sequence and
$f:=\sum_{\xi\in \cX}h_\xi\widetilde{\psi}_\xi$.
Similarly as in the proof of Theorem~\ref{thm:Fnorm-equivalence},
we use Lemmas~\ref{lem:B-maximal} and \ref{lem:Phi*psi}, and
(\ref{eq:useful}) to obtain
\begin{align*}
\WW(2^j; x)^{-\r/d}|\Phi_j\ast f(x)|
&\le c\sum_{j-1\le \nu\le j+1}\sum_{\xi\in \cX_\nu}
|h_\xi|\WW(2^j; x)^{-\r/d}|\Phi_j\ast\wt\psi_{\xi}(x)|\\
&\le c\sum_{j-1\le \nu\le j+1}\sum_{\xi\in \cX_\nu}
|h_\xi| \frac{2^{jd/2} \WW(2^j; x)^{-\r/d}}{\sqrt{ \WW(2^j;x)}(1+2^jd(\xi, x))^k } \\
&\le c2^{jd/2}\sum_{j-1\le \nu\le j+1}\sum_{\xi\in \cX_\nu}
|h_\xi| \frac{\WW(2^j; \xi)^{-\r/d-1/2}}{(1+2^jd(\xi, x))^{k-2\mu|\r|/d-\mu} } \\
&\le c2^{jd/2}\sum_{j-1\le \nu\le j+1}\sum_{\xi\in \cX_\nu} |h_\xi|
\WW(2^j; \xi)^{-\r/d-1/2}\cM_t(\ONE_{R_\xi})(x),
\end{align*}
where $\cX_{-1}:=\emptyset$ and
in the fourth inequality we used that
$k \ge 2\mu|\r|/d+\mu+(2\mu+d)/t$.
Now employing the maximal inequality (\ref{max-ineq}) we get
\begin{equation*}
\begin{aligned}
&\norm{\WW(2^j; \cdot)^{-\r/d}\Phi_j*f(\cdot)}_{\Lp}\\
&\qquad\qquad\qquad \le c2^{jd/2}\Big\|\sum_{j-1\le \nu\le j+1}\sum_{\xi\in \cX_\nu}
|h_\xi|\WW(2^j; \xi)^{-\r/d-1/2}\cM_t(\ONE_{R_\xi})(\cdot)\Big\|_{\Lp}\\
&\qquad\qquad\qquad \le c2^{jd/2}\Big\|\sum_{j-1\le \nu\le j+1}\sum_{\xi\in \cX_\nu}
|h_\xi|\WW(2^j; \xi)^{-\r/d-1/2}\ONE_{R_\xi}(\cdot)\Big\|_{\Lp}\\
%
%
&\qquad\qquad\qquad \le c2^{jd(1/2-1/p)}\Big(\sum_{j-1\le \nu\le j+1}\sum_{\xi\in \cX_\nu}
|h_\xi|^p\WW(2^j; \xi)^{-(\r/d-1/p+1/2)p}\Big)^{1/p}.
\end{aligned}
\end{equation*}
Using this in Definition~\ref{def-Besov-sp} we obtain
$\norm{f}_{\Bsrpq}\le c\norm{\{h_\xi\}}_{\bsrpq}$.

Further, we extend this result to an arbitrary sequence $h=\{h_\xi\}\in \bsrpq$
similarly as in the Triebel-Lizorkin case
by using the embedding of $\Bsrpq$ in $\cD'$
and the density of finitely supported sequences in $\bsrpq$.

\smallskip

We next prove the boundedness of the operator
$S_\varphi:\Bsrpq\rightarrow \bsrpq$,
assuming that the space $\Bsrpq$ is defined in terms
of $\{\overline{\Phi}_j\}$ in place of $\{\Phi_j\}$.
Observe first that
$$
|\ip{f,\varphi_\xi}|
\sim m(R_\xi)^{1/2}|\overline{\Phi}_j*f(\xi)|
\sim 2^{-jd/2}\WW(2^j; \xi)^{1/2}|\overline{\Phi}_j*f(\xi)|,
\quad \xi \in \cX_j.
$$
Since $\overline{\Phi}_j\ast f\in \Pi_{2^j}$, Lemma~\ref{l:half_shannon} yields
\begin{align*}
&\sum_{\xi\in \cX_j}\Big(\WW(2^j;\xi)^{-\r/d+1/p-1/2}|\ip{f,\varphi_\xi}|\Big)^p\\
&\qquad\qquad \le c2^{-jd(p/2-1)}
\sum_{\xi\in \cX_j}\WW(2^j;\xi)^{-\r p/d}|\overline{\Phi}_j*f(\xi)|^p m(R_\xi)\\
&\qquad\qquad \le c2^{-jd(p/2-1)}\norm{\WW(2^j;\xi)^{-\r/d}\overline{\Phi}_j\ast f}_{\Lp}^p.
\end{align*}
This at once yields
$
\|\{\langle f, \ph\rangle\}\|_{\bsrpq} \le c\norm{f}_{\Bsrpq}.
$

The identity $T_\psi\circ S_\varphi=Id$ follows by Proposition~\ref{prop:needlet-rep}.

The independence of $\Bsrpq$ of the particular selection of $\ha$
in the definition of  $\{\Phi_j\}$ follows from above exactly as
in the Triebel-Lizorkin case (see the proof of
Theorem~\ref{thm:Fnorm-equivalence}). $\qed$

The parameter $\r$ in the definition of the Besov spaces $\Bsrpq$ allow to
consider different scales of spaces.
A ``classical" choice of $\r$ would be $\r=0$.
However, we maintain that most natural are the spaces $\B sspq$ ($\r=s$).
The main advantages of the spaces
$\B sspq$  over $\B s\rho pq$ with $\rho\ne s$ are
that, first, they scale (are embedded) ``correctly" with respect
to the smoothness index $s$, and secondly, the right smoothness
spaces in nonlinear n-term weighted approximation from needles are
defined in terms of spaces $\B sspq$ (see \S\ref{Nonlin-app}
below).

%
\begin{proposition}\label{B-embedding}
Let $0<p\le p_1<\infty$, $0<q\le q_1\le \infty$, and
$-\infty<s_1\le s<\infty$. Then we have the continuous embedding
\begin{equation}\label{B-embed}
\B sspq \subset \B {s_1}{s_1}{p_1}{q_1} \quad\mbox{if}\quad
s/d-1/p=s_1/d-1/p_1.
\end{equation}
\end{proposition}
This embedding result follows immediately by applying inequality
(\ref{norm-relation2}).

\medskip


We finally want to link the weighted Besov spaces
to best polynomial approximation in $\Lpp$.
As in (\ref{def-En}), let $E_n(f)_p$ denote the best approximation of $f \in \Lpp$
from $\Pi_n$.


\begin{proposition}\label{prop:character-Besov}
Let $s > 0$, $1 \le p \le \infty$, and $0 < q \le \infty$.
Then $f \in \B s0pq$ if and only if
\begin{equation}\label{character-Besov1}
\|f\|_{B^{s 0}_{pq}}^A
:= \|f\|_\Lp +
\Big(\sum_{j=0}^\infty (2^{s j}E_{2^j}(f)_p)^q\Big)^{1/q}
< \infty.
\end{equation}
Moreover,
\begin{equation}\label{character-Besov2}
\|f\|_{\B s0pq}^A \sim \|f\|_{\B s0pq}.
\end{equation}
\end{proposition}

The proof of this proposition is similar to the proof of
Proposition 5.3 in \cite{NPW2} and Proposition 6.2 in \cite{KPX}.
We omit it.

\section{Application of weighted Besov spaces to nonlinear approximation}\label{Nonlin-app}
\setcounter{equation}{0}

Let us consider nonlinear n-term approximation for a needlet system
$\{\psi_\eta\}_{\eta\in \cX}$ defined as in (\ref{def.Phi-j})-(\ref{def-needlets2})
with $\hb = \ha$, $\ha\ge 0$. Thus $\ph_\eta=\psi_\eta$ are real-valued.
Then by Proposition~\ref{prop:needlet-rep}, for any $f\in\Lpp$, $1\le p\le \infty$,
$$
f=\sum_{\xi\in\cX} \langle f, \psi_\xi\rangle \psi_\xi
\quad \mbox{in } \Lpp.
$$

Suppose $\Sigma_n$ is the nonlinear set of all functions $g$ of the form
$$
g = \sum_{\xi \in \Lambda} a_\xi \psi_\xi,
$$
where $\Lambda \subset \cX$, $\#\Lambda \le n$,
and $\Lambda$ may vary with $g$.
Let $\sigma_n(f)_p$ denote the error of best $\Lpp$-approximation to
$f \in \Lpp$ from $\Sigma_n$, i.e.
$$
 \sigma_n(f)_p := \inf_{g \in \Sigma_n} \|f - g\|_p.
$$
We consider approximation in $\Lpp$, $0<p<\infty$.
Suppose $s>0$ and let $1/\tau := s/d+1/p$.
Denote briefly
$$
\Bstau:=\B ss\tau\tau.
$$
From Theorem~\ref{thm:Bnorm-eq} and \eqref{norm-Needlets} one derives the following
representation of the norm in $\Bstau$:
\begin{equation}\label{Btau-norm}
\|f\|_{\Bstau}\sim
\Big(\sum_{\xi\in\cX} \|\langle f, \psi_\xi \rangle \psi_\xi\|_p^\tau\Big)^{1/\tau}.
\end{equation}

The following embedding result shows the importance of the spaces $\Bstau$ fot
nonlinear approximation from needlets.


\begin{proposition}\label{prop:embed-Lp}
If $f \in \Bstau$, then $f$ can be identified as a function
$f\in \Lpp$ and
\begin{equation}\label{embedding}
\|f\|_p \le
\Big\|\sum_{\xi\in\cX}|\langle f, \psi_\xi\rangle\psi_\xi(\cdot)|\Big\|_p
\le c\|f\|_{\Bstau}.
\end{equation}
\end{proposition}
For the proof one proceeds exactly as in the proof of the embedding result
from \cite[Theorem 3.3]{KP} (see also \cite[Proposition 8.1]{KPX}).
The proof will be omitted.

We now give the main result of this section.


\begin{theorem}\label{thm:jackson} {\rm [Jackson estimate]}
If $f \in \Bstau$, then
\begin{equation}\label{jackson}
\sigma_n(f)_p \le c n^{-s}\|f\|_{\Bstau}.
\end{equation}
\end{theorem}

The proofs of this theorem can be carried out exactly as the proof of
Theorem~3.4 in \cite{KP} or \cite[Theorem~8.2]{KPX} and will be omitted.

Here the main open problem is to prove the companion to (\ref{jackson})
Bernstein estimate:
\begin{equation}\label{bernstein}
\|g\|_{\Bstau} \le c n^s \|g\|_p
\quad \hbox{for}\quad g \in \Sigma_n,
\quad 1<p<\infty.
\end{equation}
This estimate would allow to characterize the rates
of nonlinear n-term approximation in $\Lpp$ ($1<p<\infty$) from needlet systems.


\section{Proofs}\label{proofs}
\setcounter{equation}{0}


\subsection{Proofs for Sections~\ref{preliminaries}-\ref{def-needlets}}\label{proofsA} ${}$

\medskip

\noindent {\em Proof of Theorem~$\ref{thm:est-Lp-norm}$.} We shall
first establish (\ref{eq:est-Lp-norm}) for $p=2$.
From the definition of the kernels $\PP_n(x,y)$
(see (\ref{def-Pn})-(\ref{compactPn}))
it follows that
$$
\int_{B^d}\PP_n(x,y)\PP_m(x,y)\W(y)\, dy=\delta_{n,m}\PP_n(x,x)
$$
and hence
\begin{equation}\label{Rep-Ln}
 \int_{B^d} |L_n(x,y)|^2 \W(y)dy  = \sum_{k=0}^{2n}
   \Big|\wh a\Big(\frac{k}{n}\Big)\Big|^2 \PP_k(x, x).
\end{equation}
Therefore, for $p=2$ estimate (\ref{eq:est-Lp-norm}) will follow by
the following lemma.


\begin{lemma}\label{L2-lower-bound}
For any $\eps>0$
\begin{equation} \label{lowerbd}
\sum_{j=n}^{n+[\eps dn]} \PP_j(x,x)
\ge \frac{cn^{d}}{\WW(n;x)},
\qquad x \in B^d, \quad n\ge 1/\eps,
\end{equation}
where $c>0$ depends only on $\eps$, $\mu$, and $d$.
\end{lemma}

\begin{proof}Assume $\mu > 0$. We shall utilize
representation (\ref{compactPn}) of $\PP_n(x, y)$. The case
$\mu=0$ is easier and will be omitted (in this case one uses
representation (4.2) of $\PP_n(x, y)$ from \cite{PX2}).

From \eqref{compactPn} it is obvious that $\PP_n(x,x)$ depends only on $|x|$.
For the rest of the proof, we denote
$\PP_{n,d}(r) := \PP_n(x,x)$, where $r:=|x|$, and
$
\Lambda_{n,d}(r):= \sum_{j=n}^{n+[\eps dn]} \PP_{j, d}(r).
$
Summing up the well known recurrence relation \cite[(4.7.29)]{Sz}
$$
C_n^\lambda(x)-C_{n-2}^\lambda(x)=\frac{n+\lambda-1}{\lambda-1}C_n^{\lambda-1}(x),
\quad \hbox{where}\quad
C_{-1}^\lambda(x)=C_{-2}^\lambda(x):=0,
$$
we get
$$
C_n^\lambda(x)
=\sum_{0\le 2j\le n}\frac{n-2j+\lambda-1}{\lambda-1}C_{n-2j}^{\lambda-1}(x).
$$
Combining this with \eqref{compactPn} we arrive at
$$
\PP_{n, d}(r)= \frac{b_d^\mu}{b_{d-2}^\mu}\frac{n+\lambda}{\lambda}
\sum_{0\le 2j\le n}\PP_{n-2j, d-2}(r).
$$
Hence
\begin{align*}
\Lambda_{n,d}(r)
&=  \sum_{k=n}^{n+[\eps dn]} \PP_{k,d}(r)
= \frac{b_d^\mu}{b_{d-2}^\mu}
\sum_{k=n}^{n+[\eps dn]}
\frac{k+\lambda}{\lambda} \sum_{0 \le 2j \le k} \PP_{k-2j,d-2}(r) \\
& \ge  c \, n^2 \sum_{k=n}^{n+[\eps (d-2)n]}
\PP_{k,d-2}(r) = c \, n^2 \Lambda_{n, d-2}(r).
\end{align*}
Here $c>0$ depends only on $\eps$, $\mu$, and $d$; we used that $n\ge 1/\eps$.

Evidently, the above estimate leads to  \eqref{lowerbd} using induction on $d$,
provided we prove \eqref{lowerbd} for $d =1$ and $d=2$.
However, the case $d =1$ is already established in \cite[Proposition 2.4]{KPX},
namely,
\begin{equation}\label{est-Lambda1}
\Lambda_{n, 1}(r)\ge \frac{cn}{\WW(n;r)}.
\end{equation}

It remains to prove \eqref{lowerbd} in the case $d =2$.
The proof relies on the well known identity \cite[p. 59]{Askey}
\begin{equation}\label{rel-Gegen}
C_n^\lambda(x) = \sum_{0 \le 2k \le n}
\frac{\Gamma(\mu)(n-2k+\mu)\Gamma(k+\lambda-\mu)\Gamma(n-k+\lambda)}
     {\Gamma(\lambda)\Gamma(\lambda-\mu)k!\Gamma(n-k+\mu+1)}\,
     C_{n-2k}^\mu(x)
\end{equation}
and the product formula of Gegenbauer polynomials
\cite[Vol I, Sec. 3.15.1, (20)]{Edelyi}:
\begin{equation}\label{prod-Gegen}
 \frac{C_n^\mu(s) C_n^\mu(t)} {C_n^\mu (1)}
 = b_1^{\mu-1/2}\int_{-1}^1C_n^\mu\left(s t + u\sqrt{1-s^2}\sqrt{1-t^2}\right)
           (1-u^2)^{\mu-1} du.
\end{equation}
Using
\eqref{rel-Gegen}  (with $\lambda=\mu+1/2$)
along with \eqref{compactPn} and then \eqref{prod-Gegen}, we obtain
\begin{align}\label{est-Pn2}
\PP_{n,2}(r)
&= b_2^\mu \frac{n+\mu +1/2}{\mu+1/2} \sum_{0 \le 2k \le n}
c_{k,n} \frac{n-2k+\mu}{\mu} \frac{[C_{n-2k}^\mu(r)]^2}{C_{n-2k}^\mu(1)} \notag\\
&= \frac{b_2^\mu}{b_1^\mu} \frac{n+\mu +1/2}{\mu+1/2} \sum_{0 \le 2k \le n}
c_{k,n} \PP_{n-2k,1}(r),
\end{align}
where
$$
c_{k,n} = \frac{\Gamma(\mu+1)\Gamma(k+1/2)\Gamma(n-k + \mu +1/2)}
               {\Gamma(\mu+1/2)\Gamma(1/2)\Gamma(n-k + \mu +1) k!}.
$$
Here we used that the $L_2(\W)$-normalized Gegenbauer polynomial
$\widetilde{C}_n^\mu$ can be written in the form
$\widetilde{C}_n^\mu(x)=h_n^{-1/2}C_n^\mu(x)$ with
$h_n:= (b_1^{\mu})^{-1}\frac{\mu}{n+\mu}C_n^\mu(1)$,
which is a matter of simple verification,
and hence
$$
\PP_{n, 1}(r)
=[\widetilde{C}_n^\mu(r)]^2
=b_1^{\mu}\frac{n+\mu}{\mu}\frac{[C_n^\mu(r)]^2}{C_n^\mu(1)}.
$$
It is straightforward to verify that if $0 \le k \le n/2$,
then $c_{k,n} \sim (kn)^{-1/2}$
and hence $c_{k,n} \ge c n^{-1}$.
Therefore, from \eqref{est-Pn2}
\begin{align*}
\Lambda_{n,2}(r)
&= \sum_{k=n}^{n+[2\eps n]} \PP_{k,2}(r)
= \frac{b_2^\mu}{b_1^\mu}\sum_{k=n}^{n+[2\eps n]}\frac{k+\mu +1/2}{\mu+1/2}
\sum_{0 \le 2j \le k}c_{j,k} \PP_{k-2j,1}(r) \\
& \ge c  \sum_{k=n}^{n+[2\eps n]}
\sum_{0 \le 2j \le k}  \PP_{k-2j,1}(r)
\ge  c\, n  \Lambda_{n,1}(r).
\end{align*}
This combined with \eqref{est-Lambda1} yields \eqref{lowerbd} for
$d=2$.
\end{proof}

\medskip

We now continue with the proof of Theorem~\ref{thm:est-Lp-norm}.
Applying \eqref{lowerbd} with $\eps=2/3d$ yields
$
\|L_n(x, \cdot)\|_2 \ge cn^d\WW(n;x)^{-1}
$
for $n\ge 2d$.
If~$2\le n < 2d$, then as in the proof of Lemma~\ref{L2-lower-bound} it follows that
$$
\|L_n(x, \cdot)\|_2^{1/2} \ge c(\PP_n(x, x) + \PP_{n+1}(x,x))
\ge c(C_n^\mu(|x|) + C_{n+1}^\mu(|x|)) >c>0
$$
for all $x\in B^d$,
where we used the fact that the polynomials $C_n^\mu$ and $C_{n+1}^\mu$ have no common zeros.
Taking into account that $\WW(n;x) \sim 1$ when $n\le 2d$, the above leads again to
$
\|L_n(x, \cdot)\|_2 \ge cn^d\WW(n;x)^{-1}.
$
This completes the proof of estimate \eqref{eq:est-Lp-norm}  for $p=2$.

Now one easily derives 
\eqref{eq:est-Lp-norm} for
$p\ne 2$ from the same estimate for $p=2$ and the upper bound estimate
\eqref{eq:LpUB}.
Indeed, for $2<p<\infty$ applying H\"{o}lder's inequality we get
\begin{align*}
\frac{cn^d}{\WW(n,x)}
&\le \int_{B^d}|L_n(x, y)|^2\W(y)dy
\le \|L_n(x, \cdot)\|_p\|L_n(x, \cdot)\|_{p'}\\
&\le c_1\|L_n(x, \cdot)\|_p\Big(\frac{n^d}{\WW(n,x)}\Big)^{1-1/p'}
\quad (1/p+1/p'=1),
\end{align*}
which implies \eqref{eq:est-Lp-norm}.
One proceeds similarly whenever $p=\infty$.

If $0<p<2$, using \eqref{eq:est-Lp-norm} for $p=2$ and
\eqref{eq:LpUB} for $p=\infty$, we get
\begin{align*}
\frac{cn^d}{\WW(n,x)}
&\le \int_{B^d}|L_n(x, y)|^2\W(y)dy
\le \int_{B^d}|L_n(x, y)|^p\W(y)dy\|L_n(x, \cdot)\|_\infty^{2-p}\\
&\le c_1\int_{B^d}|L_n(x, y)|^p\W(y)dy\Big(\frac{n^d}{\WW(n,x)}\Big)^{2-p},
\end{align*}
This again leads to \eqref{eq:est-Lp-norm}. The proof of
Theorem~\ref{thm:est-Lp-norm} is complete. $\qed$

\medskip


\noindent {\em Proof of Proposition~$\ref{Nikolski}$.}
Let $g\in\Pi_n$.
Assume $1<q<\infty$ and let $L_n$ be the kernel from $(\ref{def-Ln})$,
with $\ha$ admissible of type $(a)$.
By Lemma~\ref{lem:Ker-n}, $g=L_n*g$.
We~use this, H\"older's inequality,
\eqref{eq:LpUB}, and that $\WW(n;x)\ge n^{-2\mu}$
to obtain
$$
   |g(x)| \le \|g\|_q\left(\frac{n^d}{\WW(n;x)}\right)^{1/q}
         \le c n^{(d+2 \mu) /q}\|g\|_q,
\quad x\in B^d,
$$
and hence
\begin{equation}\label{Nik1}
\|g\|_\infty \le c n^{(d+2\mu)/q} \|g\|_q,
\quad 1<q \le \infty.
\end{equation}
Let $0<q\le 1$. The above inequality with $q=2$ yields
$$
\|g\|_\infty^2
\le c n^{d+2\mu} \int_{-1}^1 |g(y)|^{2-q} |g(y)|^q \W(y)dy
   \le c n^{d+ 2\mu}\|g\|_\infty^{2-q}\|g\|_q^q.
$$
Therefore, (\ref{Nik1}) holds for $0<q\le 1$ as well.

Let $0<q<p<\infty$. Using (\ref{Nik1}) we have
\begin{align*}
\|g\|_p  & = \left(\int_{B^d} |g(x)|^{p-q} |g(x)|^q \W(x) dx \right)^{1/p}\\
   & \le c n^{(d+2\mu)(\frac{1}{q} - \frac{1}{p})} \|g\|_q^{\frac{p-q}{p} }
     \|g\|_q^{\frac{q}{p}}
    =  c n^{(d+2\mu )(\frac{1}{q} - \frac{1}{p})} \|g\|_q.
\end{align*}
Thus we have proved \eqref{norm-relation}.
$\qed$


We next prove (\ref{norm-relation2}). Assume first that $1<q<\infty$.
Using again that $g=L_n*g$, H\"older's inequality ($1/q + 1/q'=1$),
and \eqref{eq:main-estI} we obtain for $x\in B^d$,
\begin{align*}
|g(x)|
& \le \|\WW(n;\cdot)^{\gamma +\frac{1}{p} - \frac{1}{q}} g(\cdot)\|_q
\left( \int_{B^d} \Big|L_n(x,y) \WW(n;y)^{-\gamma-\frac{1}{p}+
\frac{1}{q}}\Big|^{q'} \W(y) dy \right)^{1/q'} \\
& \le c \frac{n^d} {\WW(n;x)^{1/2}} \left (\int_{B^d}
\frac{\W(y)dy}{ \WW(n;y)^{\frac{q'}{2} + \beta} (1+ n d(x,y))^\sigma } \right)^{1/q'}
\|\WW(n;\cdot)^{\gamma +\frac{1}{p} - \frac{1}{q}} g(\cdot)\|_q,
\end{align*}
where $\beta = q'(\gamma + \frac{1}{p}-\frac{1}{q})$.
The last integral can be estimated by using \eqref{eq:LpUB}, yielding
$$
 |g(x)|  \le  c \frac{n^{d/q}}{\WW(n;x)^{\gamma +1/p}}
     \|\WW(n;\cdot)^{\gamma +\frac{1}{p} - \frac{1}{q}} g(\cdot)\|_q.
$$
Hence
\begin{equation} \label{aaa2}
\|\WW(n;\cdot)^{\gamma+1/p}g(\cdot)\|_\infty \le c n^{d/q}
\|\WW(n;\cdot)^{\gamma +\frac{1}{p} - \frac{1}{q}} g(\cdot)\|_q,
\quad 1<q\le \infty.
\end{equation}
Let $0 < q \le 1$. Then by (\ref{aaa2}) with $q=2$ we have
\begin{align*}
 \|\WW(n;\cdot)^{\gamma+1/p}g(\cdot)\|_\infty
&\le c n^{d/2}
 \|\WW(n;\cdot)^{\gamma +\frac{1}{p} - \frac{1}{2}} g(\cdot)\|_2 \\
 & \le c n^{d/2}  \| \WW(n;\cdot)^{\gamma+1/p}g(\cdot)\|_\infty^{1-q/2}
 \|\WW(n;\cdot)^{\gamma +\frac{1}{p} - \frac{1}{q}} g(\cdot)\|_q^{q/2}.
\end{align*}
Therefore, \eqref{aaa2} holds for $0 < q \le 1$ as well.

Let $p<\infty$.
Using (\ref{aaa2}), we have
\begin{align*}
\|\WW & (n;\cdot)^\gamma g(\cdot)\|_p
= \left( \int_{B^d} \left[\WW(n;x)^\gamma g(x)\right]^{p-q}
\left[\WW(n;x)^\gamma g(x)\right]^{q}\W(x) dx \right)^{1/p} \\
\le & c n^{d(\frac{1}{q} - \frac{1}{p})}
\|\WW(n;\cdot)^{\gamma +\frac{1}{p} - \frac{1}{q}} g(\cdot)\|_q^{p-q}
\left( \int_{B^d} \frac{\left[\WW(n;x)^\gamma g(x)\right]^{q}}
{\WW(n;x)^{\frac{p-q}{p}}} \W(x) dx \right)^{1/p} \\
= & c n^{d(\frac{1}{q}-\frac{1}{p})}
\|\WW(n;\cdot)^{\gamma +\frac{1}{p} - \frac{1}{q}} g(\cdot)\|_q.
\end{align*}
Hence \eqref{norm-relation2} holds for $p < \infty$. If $p =
\infty$ \eqref{norm-relation2} follows from \eqref{aaa2}.
 $\qed$

\bigskip


\noindent {\em Proof of $(\ref{norm-Needlets2})$.} From
\eqref{local-Needlets21} with $k$ sufficiently large ($k >d+2\mu$
will do), and (\ref{norm-Needlets}), we infer for $0<r\le\pi$
\begin{equation*}
\begin{aligned}
0&<c_1\le \norm{\ph_\xi}_2\\
&\le \norm{\ph_\xi}_{L_\infty(B_\xi(r))}m(B_\xi(r))
+c 2^{jd}\int_{B^d\setminus B_\xi(r)}\frac{\W(y)}{\WW(2^j;y)(1+2^jd(\xi,y))^{2k}}\, dy\\
&\le \norm{\ph_\xi}_{L_\infty(B_\xi(r))}m(B_\xi(r))
+c\frac{2^{jd}}{(1+2^j r)^k}\int_{B^d}\frac{\W(y)}{\WW(2^j;y)(1+2^jd(\xi,y))^k}\, dy\\
&\le \norm{\ph_\xi}_{L_\infty(B_\xi(r))} m(B_\xi(r))+ \frac{c_2}{1+2^j r},
\end{aligned}
\end{equation*}
where $c_2$ depends only on $k$, $d$, and $\mu$.
For the last inequality we used Lemma~\ref{lem:instrumental} with $p=2$.
Let $r:=c^* 2^{-j}$, where $c^*>0$ is selected so that
$ c_2/(1+2^j r) = c_2/(1+c^*) < c_1/2$. Then from above
$$
\norm{\ph_\xi}_{L_\infty (B_\xi(c^*2^{-j}))}
\ge \frac{c}{m(B_\xi(c^*2^{-j}))}\ge c \Big(\frac{2^{jd}}{\WW(2^j;\xi)}\Big)^{1/2}.
$$
A similar estimate holds for $\psi_\xi$ as well.  $\qed$

\subsection{Proofs for Sections 4-5}\label{proofsB}${}$

\medskip

\noindent {\em Proof of Lemma $\ref{lem:Phi*psi}$.}
Using the orthogonality of the subspaces $\cV_n^d$, we have
$\Phi_j\ast \psi_\xi(x)=0$ if $\xi\in\cX_\nu$ for $\nu\ge j+2$ or
$\nu\le j-2$.

Let $\xi \in \cX_\nu$, $j-1\le \nu\le j+1$.
From the localization of the kernels $\Phi_j$, given in (\ref{local-Needlets2}),
and the needlet localization from (\ref{local-Needlets21}) it follows that for any $k>0$
there is a constant $c_k>0$ such that
\begin{equation*}
\begin{aligned}
&|\Phi_j\ast \psi_\xi(x)|
\le  c_k \frac{2^{j3d/2}}{\sqrt{\WW(2^j;x)}}\int_{B^d}\frac{\W(y)}
{\sqrt{ \WW(2^j;y)}(1+2^jd(x,y))^k(1+2^jd(y,\xi))^k}\, dy.
\end{aligned}
\end{equation*}
Denote
$$
\Omega_\xi:=\{y\in B^d: d(y,\xi)\ge d(x,\xi)/2\} \text{ and }
\Omega_x:=\{y\in B^d: d(x, y)\ge d(x,\xi)/2\}.
$$
Evidently, $B^d=\Omega_{\xi}\cup \Omega_x$ and hence
\begin{align*}
|\Phi_j\ast \psi_\xi(x)|
&\le c_k\frac{2^{j3d/2}}{\sqrt{ \WW(2^j;x)}(1+2^jd(x, \xi))^k} \int_{\Omega_\xi}
 \frac{\W(y)}{\WW(2^j; y)(1+2^jd(x,y))^k}\, dy\\
&+ c_k\frac{2^{j3d/2}}{\sqrt{ \WW(2^j;x) }(1+2^jd(x, \xi))^k}
\int_{\Omega_x} \frac{ \W(y)}{\WW(2^j; y)(1+2^jd(y,\xi))^k}\, dy\\
&=:J_1+J_2.
\end{align*}
We may assume that $k>d$. Then employing Lemma~\ref{lem:instrumental} with $p=2$,
we get
\begin{align*}
\int_{\Omega_\xi}
\frac{ \W(y)}{\WW(2^j; y)(1+2^jd(x,y))^k}\, dy\le \int_{B^d}
\frac{ \W(y)}{\WW(2^j; y)(1+2^jd(x,y))^k}\, dy\le c2^{-jd},
\end{align*}
which yields
$$
J_1\le c\frac{2^{jd/2}}{\sqrt{\WW(2^j;x)}(1+2^jd(x, \xi))^k}.
$$
One similarly estimates $J_2$. This completes the proof of the
lemma.  $\qed$

\medskip


\noindent {\em Proof of Lemma~$\ref{lem:Max-needl}$.}
 Estimate \eqref{Max-needl1} follows readily from the localization of
 the needlets (see (\ref{local-Needlets21})) and the lower bound estimate from (\ref{B-max3})
 taking into account that
 $R_\xi \subset B_\xi(c^\diamond 2^{-j})$ for $\xi\in\cX_j$.

We now prove (\ref{Max-needl2}).
By the lower bound estimate (\ref{norm-Needlets2}) it follows that there exists
$\omega\in B_\xi(c^*2^{-j})$ such that
\begin{equation}\label{eq:needlUB}
|\ph_\xi(\omega)|\ge c  \frac{2^{jd/2}}{\sqrt{\WW(2^j; \xi)}}.
\end{equation}
Also, by (\ref{Lip-Needlets}) it follows that for every $x\in B_\omega(2^{-j})$
\begin{equation}\label{eq:needl-smooth}
 |\ph_\xi(\omega)-\ph_\xi(x)|
 \le c\frac{2^{j(d/2+1)}d(\omega,x)}{\sqrt{\WW(2^j; \xi)}}.
\end{equation}
By (\ref{eq:needlUB})-(\ref{eq:needl-smooth}) it follows that
for a sufficiently small constant $\hat{c}>0$
$$
|\ph_\xi(x)|
\ge |\ph_\xi(\omega)|-|\ph_\xi(\omega)-\ph_\xi(x)|
\ge c\frac{2^{jd/2}}{\sqrt{\WW(2^j; \xi)}}
\ge c\tONE_{B_\omega(\hat{c}2^{-j})}(x),
\quad x\in B_\omega(\hat{c}2^{-j}),
$$
which yields
\begin{equation*}
(\cM_t\ph_\xi)(x)
\ge c (\cM_t\tONE_{B_\omega(\hat{c}2^{-j})})(x)
\ge c\tONE_{B_\xi(2^{-j})}(x)
\ge c\tONE_{R_\xi}(x),
\quad x\in B^d,
\end{equation*}
where in the second inequality we used (\ref{B-max3}).

One similarly shows that $\cM_t\psi_\xi \ge c\tONE_{R_\xi}$.
 $\qed$

\medskip


\noindent {\em Proof of Lemma~$\ref{lem:weak_inequality}$.} For
the proof of this lemma we need a couple of additional lemmas.


\begin{lemma}\label{lem:sums}
Let $k>d$ and $j\ge 0$. Then
\begin{equation} \label{eq:sum1}
\sum_{\xi\in \cX_j} \frac{1}{(1+2^jd(x,\xi))^k}\le c,
\quad x\in B^d,
\end{equation}
and for any $\xi,\eta\in B^d$
\begin{equation}\label{eq:sum2}
\sum_{w\in \cX_j}\frac{1}{(1+2^jd(\xi,w))^k(1+2^jd(\eta,w))^{k}}
\le c\frac{1}{(1+2^jd(\xi,\eta))^{k}}.
\end{equation}
\end{lemma}

\begin{proof}Fix $\xi\in\cX_j$. Evidently, $1+2^jd(x,
\xi)\sim 1+2^jd(x, y)$ for $y \in R_\xi$, and by
\eqref{norm-dist2}
$$
|\sqrt{1-|\xi|^2}-\sqrt{1-|y|^2}|\le \sqrt{2}\,d(\xi, y)\le c2^{-j},
\quad y\in R_\xi,
$$
which implies
$$
|R_\xi|\sim 2^{-jd}(\sqrt{1-|\xi|^2}+2^{-j})
\sim 2^{-jd}(\sqrt{1-|y|^2}+2^{-j}),
\quad y\in R_\xi.
$$
We use the above to obtain
\begin{align*}
\sum_{\xi\in\cX_j}\frac{1}{(1+2^jd(x, \xi))^k}
&\le c\sum_{\xi\in\cX_j}\frac{1}{|R_\xi|} \int_{R_\xi}\frac{1}{(1+2^jd(x,y))^k}\,dy\\
&\le c2^{jd}\int_{B^d}\frac{1}{(\sqrt{1-|y|^2}+2^{-j})(1+2^jd(x,y))^k}\,dy
\le c.
\end{align*}
Here for the last inequality we used Lemma~\ref{lem:instrumental}
with $p=2$ and $\mu=1/2$.


For the proof of \eqref{eq:sum2}, assume that $\xi\ne \eta$ and denote
$$
\cX_j(\xi):=\{w\in \cX_j: d(\xi,w)\ge d(\xi,\eta)/2\}, \quad
\cX_j(\eta):=\{w\in \cX_j: d(\eta,w)\ge d(\xi,\eta)/2\}.
$$
Then
\begin{align*}
&\sum_{w\in \cX_j}\frac{1}{(1+2^jd(\xi,w))^k(1+2^jd(\eta,w))^{k}}
\le c\frac{1}{(1+2^jd(\xi,\eta))^k}
\sum_{w\in \cX_j(\xi)}\frac{1}{(1+2^jd(\eta,w))^{k}}\\
&\qquad\qquad\quad + c\frac{1}{(1+2^jd(\xi,\eta))^k}
\sum_{w\in \cX_j(\eta)}\frac{1}{(1+2^jd(\xi,w))^{k}}\\
&\qquad\qquad\quad \le c\frac{1}{(1+2^jd(\xi,\eta))^k}
\Big(\sum_{w\in \cX_j}\frac{1}{(1+2^jd(\eta,w))^{k}}
+\sum_{w\in \cX_j}\frac{1}{(1+2^jd(\xi,w))^{k}}\Big)\\
&\qquad\qquad\quad \le c\frac{1}{(1+2^jd(\xi,\eta))^k},
\end{align*}
where for the last inequality we used \eqref{eq:sum1}.
 \end{proof}


\begin{lemma}\label{lem:quasi-Lip}
Assume $P\in \Pi_{2^j}$ $(j \ge 0)$, $\xi\in \cX_j$, and
let $x_1,x_2\in B^d$ and
$d(x_\nu,\eta)\le \tilde{c}2^{-j}$, $\nu=1,2$.
For any $k>0$
$$
|P(x_1)-P(x_2)|\le c 2^j d(x_1, x_2)
\sum_{\xi\in\cX_j}\frac{|P(\xi)|}{(1+2^j d(\eta, \xi))^k},
$$
where $c>0$ depends only on $d$, $k$, $\mu$, and $\tilde{c}$.
\end{lemma}

\begin{proof}Fix $P\in\Pi_{2^j}$ and assume that
$L_{2^j}$ is the reproducing kernel from Lemma~\ref{lem:Ker-n}
with $n=2^j$. Then, $L_{2^j}*P=P$. Since $L_{2^j}(x,
\cdot)P(\cdot) \in \Pi_{2^{j+2}}$, and the cubature formula
(\ref{eq:cubature}) is exact for all polynomials from
$\Pi_{2^{j+2}}$ we have
$$
P(x)=\int_{B^d} L_{2^j}(x,y)P(y)\W(y)dy
=\sum_{\xi\in\cX_j}\lambda_\xi L_{2^j}(x,\xi)P(\xi),
\quad x\in B^d.
$$
We use \eqref{Lip} to obtain
for $x_1, x_2\in B^d$ with $d(x_\nu, \eta)\le \tilde{c}2^{-j}$, $\nu=1, 2$,
\begin{align*}
|P(x_1)-P(x_2)|
&=\Big|\int_{B^d}[\Ker_{2^j}(x_1,y)-\Ker_{2^j}(x_2,y)]P(y)\W(y) \, dy\Big|\\
&\le\sum_{\xi\in \cX_j}
|\lambda_{\xi}||\Ker_{2^j}(x_1,\xi)- \Ker_{2^j}(x_2,\xi)||P(\xi)|\\
&\le c2^{j} d(x_1,x_2) \sum_{\xi\in \cX_j}
\Big(\frac{\WW(2^j;\xi)}{\WW(2^j;\eta)}\Big)^{1/2}
\frac{|P(\eta)|}{(1+2^jd(\xi,\eta))^{k}}\\
&\le c2^j d(x_1,x_2) \sum_{\eta\in \cX_j}
\frac{|P(\eta)|}{(1+2^jd(\xi,\eta))^{k-2\mu}}.
\end{align*}
Here we used that $\lambda_\xi\sim 2^{-jd}\WW(2^j;\xi)$ and
for the last inequality we used (\ref{eq:useful}).
Taking into account that $k>0$ can be arbitrarily large the result follows.
 \end{proof}

\medskip

\noindent {\em Completion of the proof of Lemma
$\ref{lem:weak_inequality}$.} Since $b_\xi\le a_\xi,$ it trivially
follows that $b_\xi^*\le a_\xi^*$.

For the other direction let
$$
d_\xi:=\max\{|P(x_1)-P(x_2)|: x_1\in R_\xi, d(x_1,x_2)\le  2^{-j-r}\}.
$$
Obviously $a_\xi\le b_\xi +d_\xi$.
Now Lemma \ref{lem:quasi-Lip} yields
$$
d_\xi\le c2^{-r} \sum_{\eta\in\cX_j}\frac{|P(\eta)|}{(1+2^j d(\xi, \eta))^k},
\quad \xi\in \cX_j.
$$
From the definition of $d_\xi^*$ in (\ref{def.h-star}) we infer
\begin{align*}
d_\xi^\ast&\le c 2^{-r} \sum_{w\in \cX_j}\sum_{\eta\in\cX_j}
\frac{|P(\eta)|}{(1+2^j d(w,\eta))^k(1+2^jd(\xi,w))^k}\\
&\le c 2^{-r}\sum_{\eta\in \cX_j}\frac{|P(\eta)|}{(1+2^j d(\eta, \xi))^k}
\le c 2^{-r}a_\xi^*,
\end{align*}
where for the second inequality we interchanged the order of summation and
used Lemma \ref{lem:sums}.
Hence,
$a_\xi^*\le b_\xi^*+d_\xi^* \le b_\xi^*+c2^{-r}a_\xi^*$ with $c>0$
independent of $r$.
By selecting $r$ sufficiently large we get
$a_\xi^*\le cb_\xi^\ast$.
 $\qed$

\medskip


\noindent {\em Proof of Lemma $\ref{lem:disc-max}$.} We first
prove Lemma \ref{lem:disc-max} in the case $\r=0$. We fix
$\xi\in\cX_j$ and define $ S_0:=\{\eta\in \cX_j: d(\eta,\xi)\le
c^\diamond 2^{-j}\} $ and
$$
S_m :=\{\eta\in \cX_j: c^\diamond 2^{-j+m-1}
< d(\eta,\xi)\le c^\diamond 2^{-j+m}\}, \quad m\ge 1,
$$
where $c^\diamond$ is the constant from Proposition \ref{prop:cubature}.
By Definition~\ref{def.eq-dist-p} it follows that $\# S_m\le c2^{md}$.
Let us also set
$$
B_m := B_\xi(c^\diamond(2^m+1)2^{-j}), \quad m\ge 0.
$$
Evidently,
$R_\eta\subset B_m$ for $\eta\in S_\nu$, $0\le \nu\le m$.
Moreover, if $\eta\in S_m$, then
$$
d(\xi,\partial B^d)\le d(\xi,\eta)+d(\eta,\partial B^d)
\le c^\diamond 2^{-j+m}+d(\eta,\partial B^d).
$$
Hence, using \eqref{ball-size}, we get
\begin{equation}\label{quotient}
\begin{aligned}
\frac{m(B_m)}{m(R_\eta)}
&\le 2^{md}\left(\frac{d(\xi,\partial B^d)+ 2^{-j+m}}{d(\eta,\partial B^d)+ 2^{-j}}\right)^{2\mu}\\
&\le c 2^{md}\left(\frac{d(\eta,\partial B^d)+ 2^{-j+m}}{d(\eta,\partial B^d)+ 2^{-j}}\right)^{2\mu}
\le c 2^{m(d+2\mu)}.
\end{aligned}
\end{equation}

Set $\gamma:=\max \{0, 1-\frac1{t}\}<1$.
Using H\"older's inequality if $t> 1$ and
the $t$-triangle inequality if $0<t\le1$, we get
\begin{align*}
b_\xi^{\ast}=
\sum_{\eta\in\cX_j}\frac{|b_\eta|}{(1+2^jd(\eta,\xi))^\sigma}
\le c\sum_{m\ge 0}2^{-m\sigma}\sum_{\eta\in S_m} |b_\eta|
\le c\sum_{m\ge 0}2^{-m(\sigma-d\gamma)}(\sum_{\eta\in S_m} |b_\eta|^t)^{1/t}.
\end{align*}
We now use (\ref{quotient}) to obtain, for $x\in R_\xi$,
\begin{align*}
b_\xi^{\ast}
&=c \sum_{m=0}^\infty 2^{-m(\sigma-d)}\Bigl(\int_{B^d}
\Big[\sum_{\eta\in S_m}
|b_\eta|m(R_\eta)^{-1/t}\ONE_{R_\eta}(y)\Big]^t\W(y)\, dx\Bigr)^{1/t}\\
&\le c \sum_{m=0}^\infty 2^{-m(\sigma-d)}\Bigl(\frac1{m(B_m)}\int_{B_m}
\Big[\sum_{\eta\in S_m}
\Big(\frac{m(B_m)}{m(R_\eta)}\Big)^{1/t}|b_\eta|\ONE_{R_\eta}(y)\Big]^t\W(y)\,dy\Bigr)^{1/t}\\
&\le c\sum_{m\ge 0}2^{-m(\sigma-d-(d+2\mu)/t)}
\Bigl(\frac1{m(B_m)}\int_{B_m}
\Big[\sum_{\eta\in S_m}|b_\eta|\ONE_{R_\eta}(y)\Big]^t\W(y)\, dy\Bigr)^{1/t}\\
&\le c\cM_t\Big(\sum_{w\in \cX_j}|b_w|\ONE_{R_\omega}\Big)(x),
\end{align*}
where for the last inequality we used that $\sigma>d+(d+2\mu)/t$.

Consider now the general case.
Using (\ref{eq:useful}) we have for $\xi\in\cX_j$
\begin{align*}
\WW(2^j;\xi)^{\gamma}b_\xi^*
&\le \sum_{\eta\in\cX_j}\frac{\WW(2^j;\xi)^{\gamma}|b_\eta|}{(1+2^jd(\xi, \eta))^\sigma}
\le c\sum_{\eta\in\cX_j}\frac{\WW(2^j;\eta)^{\gamma}|b_\eta|}
                              {(1+2^jd(\xi, \eta))^{\sigma-2\mu |\gamma|}}\\
&\le c\Big(\WW(2^j;\xi)^{\gamma}|b_\xi|\Big)^*,
\end{align*}
where we used that $\sigma>d+(d+2\mu)/t+2\mu |\gamma|$. Now
(\ref{disc-max}) in the general case follows by the same
inequality in the case $\rho=0$ established above.  $\qed$

\medskip


\noindent {\em Proof of Lemma $\ref{l:half_shannon}$.} For any
$\xi\in\cX_j$, we denote $a_\xi:=\max_{x\in R_\xi}|P(x)|$,
$$
m_\xi:=\min_{x\in R_\xi}|P(x)|,
\quad \mbox{ and }\quad
b_\xi:=\max\{\min_{x\in R_w }|P(x)|:w\in \cX_{j+r},R_w\cap R_\xi\ne\emptyset \},
$$
where $r\ge 1$ is the constant from Lemma \ref{lem:weak_inequality}.

Choose $0<t<p$. By Lemma ~\ref{lem:weak_inequality} we have $a_\xi^* \le cb_\xi^*$.
We use this, Lemmas~\ref{lem:disc-max}, and the maximal inequality \eqref{max-ineq}
to obtain
\begin{equation}
\begin{aligned}\label{Wa<Wb}
\Big(\sum_{\xi\in \cX_j}\WW(2^j; \xi)^\gamma a_\xi^p m(R_\xi)\Big)^{1/p}
&=\Big\|\sum_{\xi\in \cX_j}\WW(2^j; \xi)^\gamma a_\xi \ONE_{R_\xi}(\cdot)\Big\|_{\Lp}\\
\le c\Big\|\sum_{\xi\in \cX_j}\WW(2^j; \xi)^\gamma b_\xi^\ast \ONE_{R_\xi}(\cdot)\Big\|_{\Lp}
&\le c\Big\|\cM_t\Big(\sum_{\xi\in \cX_j}\WW(2^j; \xi)^\gamma b_\xi  \ONE_{R_\xi}\Big)(\cdot)\Big\|_{\Lp}\\
&\le c\Big\|\sum_{\xi\in \cX_j}\WW(2^j; \xi)^\gamma b_\xi \ONE_{R_\xi}(\cdot)\Big\|_{\Lp}.
\end{aligned}
\end{equation}
Now, exactly as in the proof of Theorem~\ref{thm:Fnorm-equivalence} (see \eqref{W-b-xi-m-eta})
we have
\begin{equation}\label{eq:important}
b_\xi\WW(2^j; \xi)^\gamma\ONE_{R_\xi}
\le \sum_{\eta\in \cX_{j+r}(\xi)}  m^\ast_\eta \WW(2^{j+r}; \eta)^\gamma \ONE_{R_\eta}.
\end{equation}
where
$\cX_{j+r}(\xi):=\{w\in \cX_{j+r}:R_w\cap R_\xi\ne \emptyset\}$.
Combining this with (\ref{Wa<Wb}) and using that $\# \cX_{j+r}(\xi)\le c$,
Lemmas~\ref{lem:disc-max}, and the maximal inequality \eqref{max-ineq},
we get
\begin{equation*}
\begin{aligned}
\Big(\sum_{\xi\in \cX_j}\WW(2^j; \xi)^\gamma a_\xi^p m(R_\xi)\Big)^{1/p}
&\le c\Big\|\sum_{\eta\in \cX_{j+r}}m^\ast_\eta \WW(2^{j+r}; \eta)^\gamma\ONE_{R_\eta}(\cdot)\Big\|_\Lp\\
\le c\Big\|\cM_t\Big(\sum_{\eta\in \cX_{j+r}}m_\eta\WW(2^{j+r}; \eta)^\gamma\ONE_{R_\eta}\Big)(\cdot)\Big\|_\Lp
&\le c\Big\|\sum_{\eta\in \cX_{j+r}}m_\eta \WW(2^{j+r}; \eta)^\gamma \ONE_{R_\eta}(\cdot)\Big\|_\Lp\\
&\le c\norm{P}_\Lp.
\end{aligned}
\end{equation*}
Here for the forth inequality we used that $\WW(2^{j+r}; \eta)\sim
\WW(2^{j}; x)$ if $x\in \R_\eta$, $\eta\in\cX_{j+r}$. $\qed$

\end{document}